\documentclass[11pt]{article} % Stat tech report
\newcommand{\doctype}{TECH}

\usepackage[numbers]{natbib}
\RequirePackage[OT1]{fontenc}
\RequirePackage{amsthm,amsmath,amsfonts}
\RequirePackage[colorlinks]{hyperref}
\RequirePackage{epsfig, graphicx, color}
\RequirePackage{times}
\RequirePackage{latexsym}
\RequirePackage{psfrag}

\usepackage{fullpage}

\usepackage{ifthen}

\ifthenelse{\equal{\doctype}{TECH}}
{
\setlength{\topmargin}{0 in}
\setlength{\textwidth}{6 in}
\setlength{\textheight}{8.5 in}
\setlength{\oddsidemargin}{0.5 in}

}{}

\ifthenelse{\equal{\doctype}{AOS}}
{
\startlocaldefs
\endlocaldefs
}{}

\hypersetup{citecolor=red}

\theoremstyle{plain}

\newtheorem{theo}{Theorem}[section]

\newtheorem{lem}{Lemma}[section]
\newtheorem{prop}{Proposition}[section]
\newtheorem{cor}{Corollary}[section]

\theoremstyle{definition} 

\newtheorem{nota}{Notation}[section]
\newtheorem{de}{Definition}[section]
\newtheorem{exa}{Example}[section]
\newtheorem{as}{Assumption}[section]
\newtheorem{alg}{Algorithm}[section]

\newcommand{\btheo}{\begin{theo}}
\newcommand{\bde}{\begin{de}}
\newcommand{\ble}{\begin{lem}}
\newcommand{\bpr}{\begin{prop}}
\newcommand{\bno}{\begin{nota}}
\newcommand{\bex}{\begin{exa}}
\newcommand{\bcor}{\begin{cor}}
\newcommand{\spro}{\begin{proof}}
\newcommand{\bas}{\begin{as}}
\newcommand{\balg}{\begin{alg}}

\newcommand{\etheo}{\end{theo}}
\newcommand{\ede}{\end{de}}
\newcommand{\ele}{\end{lem}}
\newcommand{\epr}{\end{prop}}
\newcommand{\eno}{\end{nota}}
\newcommand{\eex}{\end{exa}}
\newcommand{\ecor}{\end{cor}}
\newcommand{\fpro}{\end{proof}}
\newcommand{\eas}{\end{as}}
\newcommand{\ealg}{\end{alg}}

\theoremstyle{plain}

\newtheorem{theos}{Theorem}
\newtheorem{props}{Proposition}
\newtheorem{lems}{Lemma}
\newtheorem{cors}{Corollary}

\theoremstyle{definition}
\newtheorem{exas}{Example}
\newtheorem{algs}{Algorithm}
\newtheorem{asss}{Assumption}
\newtheorem{defns}{Definition}

\newcommand{\btheos}{\begin{theos}}
\newcommand{\etheos}{\end{theos}}
\newcommand{\bprops}{\begin{props}}
\newcommand{\eprops}{\end{props}}
\newcommand{\bdes}{\begin{defns}}
\newcommand{\edes}{\end{defns}}
\newcommand{\blems}{\begin{lems}}
\newcommand{\elems}{\end{lems}}
\newcommand{\bcors}{\begin{cors}}
\newcommand{\ecors}{\end{cors}}
\newcommand{\bexs}{\begin{exas}}
\newcommand{\eexs}{\end{exas}}
\newcommand{\balgs}{\begin{algs}}
\newcommand{\ealgs}{\end{algs}}
\newcommand{\bass}{\begin{asss}}
\newcommand{\eass}{\end{asss}}

\newcommand{\y}{{\ensuremath{Y}}}
\newcommand{\w}{{\ensuremath{w}}}
\newcommand{\Rq}{{\ensuremath{R_q}}}

\newcommand{\Bq}{{\ensuremath{\mathbb{B}_q(\Rq)}}}
\newcommand{\s}{{\ensuremath{s}}}
\newcommand{\BZ}{{\ensuremath{\mathbb{B}_0(\s)}}}

\newcommand{\Cqr}{{\ensuremath{C_{q,p}}}}

\newcommand{\LqrOne}{{\ensuremath{L_{\qpar, \rpar}}}}
\newcommand{\UqrOne}{{\ensuremath{U_{\qpar,\rpar}}}}

\newcommand{\Ker}{\ker}

\newcommand{\diam}{\operatorname{diam}}
\newcommand{\entnum}{\ensuremath{e}}
\newcommand{\betahat}{\ensuremath{\widehat{\beta}}}
\newcommand{\betaRand}{\ensuremath{B}}
\newcommand{\numobs}{\ensuremath{n}}
\newcommand{\pdim}{\ensuremath{d}}
\newcommand{\defn}{\ensuremath{:\,=}}
\newcommand{\betaTr}{\ensuremath{\beta^{*}}}
\newcommand{\Ball}{\ensuremath{\mathbb{B}}}

\newcommand{\qpar}{\ensuremath{q}}
\newcommand{\rpar}{\ensuremath{p}}

\newcommand{\epskey}{\ensuremath{\epsilon_\numobs}}

\newcommand{\epspack}{\ensuremath{\delta_{\numobs}}}

\newcommand{\Xmat}{\ensuremath{X}}
\newcommand{\LTwo}{\ensuremath{\ell_2}}
\newcommand{\Ltwo}{\ensuremath{\LTwo}}

\newcommand{\Kx}{\ensuremath{\mathcal{N}_{\qpar}(\Xmat)}}
\newcommand{\KZx}{\ensuremath{\mathcal{N}_{0}(\Xmat)}}

\newcommand{\kOne}{\ensuremath{\nu}}

\newcommand{\mprob}{\ensuremath{\mathbb{P}}}

\newcommand{\betastar}{\ensuremath{\beta^*}}
\newcommand{\real}{\ensuremath{\mathbb{R}}}

\newlength{\widebarargwidth}
\newlength{\widebarargheight}
\newlength{\widebarargdepth}

\newcommand{\matsnorm}[2]{|\!|\!| #1 | \! | \!|_{{#2}}}

\newcommand{\betatil}{\ensuremath{\widetilde{\beta}}}
\newcommand{\Ballq}{\ensuremath{\mathbb{B}_\qpar}}
\newcommand{\myrad}{\ensuremath{R_\qpar}}

\newcommand{\mykernel}{\ensuremath{\operatorname{Ker}}}

\newcommand{\conlower}{\ensuremath{\kappa_\ell}}
\newcommand{\conupper}{\ensuremath{\kappa_u}}

\newcommand{\conlowersq}{\ensuremath{\kappa^2_\ell}}
\newcommand{\conuppersq}{\ensuremath{\kappa^2_u}}

\newcommand{\SpecFunLowerZ}{\ensuremath{f_\ell(\s, \numobs, \pdim)}}
\newcommand{\SpecFunLower}{\ensuremath{f_\ell(\Rq, \numobs, \pdim)}}
\newcommand{\qRate}{\ensuremath{g(\Rq, \numobs, \pdim)}}
\newcommand{\qRatesq}{\ensuremath{g^2(\Rq, \numobs, \pdim)}}

\newcommand{\colnorm}{\ensuremath{\kappa_c}}
\newcommand{\colnormsq}{\ensuremath{\kappa^2_c}}

\newcommand{\Annoy}{\ensuremath{c}}
\newcommand{\AnnoyLower}{\ensuremath{\Annoy_{\qpar, \rpar}}}
\newcommand{\AnnoyUpper}{\ensuremath{\Annoy'_{\qpar, 2}}}
\newcommand{\AnnoyLowerZero}{\ensuremath{\Annoy_{0, \rpar}}}

\newcommand{\AnnoyPred}{\ensuremath{\widetilde{\Annoy}_{\qpar,2}}}

\newcommand{\Prob}{\ensuremath{\mathbb{P}}}
\newcommand{\Exs}{\ensuremath{\mathbb{E}}}

\newcommand{\spindex}{\ensuremath{\s}}

\newcommand{\PackNum}{\ensuremath{M}}
\newcommand{\CovNum}{\ensuremath{N}}
\newcommand{\kull}[2]{\ensuremath{D(#1 \, \| \, #2)}}

\newcommand{\MyLqTwo}{\ensuremath{L_{\qpar,2}}}
\newcommand{\MyUq}{\ensuremath{U_{\qpar,\rpar}}}

\newcommand{\Hyper}{\ensuremath{\mathcal{H}}}
\newcommand{\HyperSub}{\ensuremath{\widetilde{\Hyper}}}
\newcommand{\Ind}{\ensuremath{\mathbb{I}}}

\newcommand{\Sset}{\ensuremath{S}}

\newcommand{\order}{\ensuremath{\mathcal{O}}}

\newcommand{\DeltaHat}{\widehat{\Delta}}

\newcommand{\GenSet}{\ensuremath{\mathcal{S}}}

\newcommand{\mytheta}{\ensuremath{\theta}}

\newcommand{\wnoise}{\ensuremath{w}}

\newcommand{\newrad}{\ensuremath{r}}
\newcommand{\Event}{\ensuremath{\mathcal{E}}}
\newcommand{\radtwo}{\ensuremath{r}}

\newcommand{\InterSet}{\ensuremath{\mathbb{S}}}
\newcommand{\gencon}{\ensuremath{c}}

\newcommand{\gnoise}{\ensuremath{\varepsilon}}

\newcommand{\thetastar}{\ensuremath{\theta^*}}
\newcommand{\ThetaSet}{\ensuremath{\Theta}}

\newcommand{\vtil}{\ensuremath{\tilde{v}}}

\newcommand{\CovMat}{\ensuremath{\Sigma}}
\newcommand{\CovMatSqrt}{\ensuremath{\CovMat^{1/2}}}

\newcommand{\rhomax}{\ensuremath{\rho}}

\newcommand{\vbre}{\ensuremath{\breve{v}}}

\newcommand{\mycon}{\ensuremath{\mu}}

\newcommand{\tupper}{\ensuremath{t_u}}
\newcommand{\tlower}{\ensuremath{t_\ell}}

\newcommand{\Tail}{\ensuremath{\mathcal{T}}}

\newcommand{\plaincon}{\ensuremath{c}}

\newcommand{\qprob}{\ensuremath{\mathbb{Q}}}

\long\def\comment#1{}

\newcommand{\rhohamm}{\ensuremath{\rho_H}}

\newcommand{\dyaeps}{\ensuremath{\epsilon}}

\newcommand{\absconv}{\ensuremath{\operatorname{absconv}}}

\newcommand{\Gue}{Gu\'{e}don }

\newcommand{\myparagraph}[1]{\paragraph{#1:}}

\newcommand{\kdim}{\ensuremath{k}}

\newcommand{\wtil}{\ensuremath{\widetilde{w}}}
\newcommand{\Xtil}{\ensuremath{\widetilde{\Xmat}}}
\newcommand{\colnormtil}{\ensuremath{\widetilde{\colnorm}}}
\newcommand{\colnormtilsq}{\ensuremath{\widetilde{\colnorm}^2}}

\newcommand{\Loss}{\ensuremath{\mathcal{L}}}

\newcommand{\Risk}{\ensuremath{\mathcal{R}}}

\newcommand{\Aset}{\ensuremath{\mathcal{A}}}

\newcommand{\BallZ}{\ensuremath{\Ball_0}}
\begin{document}

\ifthenelse{\equal{\doctype}{AOS}}{
	\begin{frontmatter}

	\title{ Minimax rates of estimation for high-dimensional
		linear regression over
		\lowercase{$\ell_q$}-balls}
	
	\begin{aug}
	  \author{Garvesh Raskutti,%\thanksref{t2}
	    \ead[label=e1]{garveshr@stat.berkeley.edu}}
	  \author{Martin J. Wainwright,
	    \ead[label=e2]{wainwrig@eecs.berkeley.edu}}
	  	\and
		  \author{Bin Yu
		    \ead[label=e3]{binyu@stat.berkeley.edu}}

\end{aug}

\begin{abstract}
Consider the standard linear regression model $\y = \Xmat \betastar +
w$, where $\y \in \real^\numobs$ is an observation vector, $\Xmat \in
\real^{\numobs \times \pdim}$ is a design matrix, $\betastar \in
\real^\pdim$ is the unknown regression vector, and $w \sim
\mathcal{N}(0, \sigma^2 I)$ is additive Gaussian noise. This paper
studies the minimax rates of convergence for estimation of $\betastar$
for $\ell_\rpar$-losses and in the $\ell_2$-prediction loss, assuming
that $\betastar$ belongs to an $\ell_{\qpar}$-ball $\Ballq(\myrad)$
for some $\qpar \in [0,1]$.  We show that under suitable regularity
conditions on the design matrix $\Xmat$, the minimax error in
$\ell_2$-loss and $\ell_2$-prediction loss scales as $\Rq \big(\frac{\log
\pdim}{n}\big)^{1-\frac{\qpar}{2}}$. In addition, we provide lower
bounds on minimax risks in $\ell_{\rpar}$-norms, for all $\rpar \in
[1, +\infty], \rpar \neq \qpar$.  Our proofs of the lower bounds are
information-theoretic in nature, based on Fano's inequality and
results on the metric entropy of the balls $\Ballq(\myrad)$, whereas
our proofs of the upper bounds are direct and constructive, involving
direct analysis of least-squares over $\ell_{\qpar}$-balls. For the
special case $q = 0$, a comparison with $\ell_2$-risks achieved by
computationally efficient $\ell_1$-relaxations reveals that although
such methods can achieve the minimax rates up to constant factors,
they require slightly stronger assumptions on the design matrix
$\Xmat$ than algorithms involving least-squares over the
$\ell_0$-ball. 
%In particular, we show that the conditions on $\Xmat$
%required by such least-squares over $\ell_{\qpar}$-ball estimators are
%satisfied with high probability by classes of non-i.i.d. Gaussian
%random matrices for which conditions such as the restricted isometry
%property (RIP), or the weaker sparse eigenvalue condition are violated
%with high probability.
\end{abstract}

\end{frontmatter}
}
{}
%%%%%%%% END AOS FRONTMATTER %%%%%%%%%%%%%%%%%%%%%%%%%%%%%%%%%%%%%%

%% BEGIN STAT TECH FRONT PAGE
\ifthenelse{\equal{\doctype}{TECH}}{

\begin{center}
{\bf{\LARGE{Minimax rates of estimation \\ for high-dimensional linear
		regression over \lowercase{$\ell_q$}-balls}}}
	
\vspace*{.3in}

\begin{tabular}{cc}
Garvesh Raskutti$^1$ &  Martin J. Wainwright$^{1,2}$ \\
	    \texttt{garveshr@stat.berkeley.edu} &
	    \texttt{wainwrig@stat.berkeley.edu} 
\end{tabular}

\vspace*{.1in}

\begin{tabular}{c}
Bin Yu$^{1,2}$ \\
\texttt{binyu@stat.berkeley.edu} 
\end{tabular}

\vspace*{.2in}

\begin{tabular}{c}
Departments of Statistics$^1$, and EECS$^2$ \\
UC Berkeley,  Berkeley, CA 94720
\end{tabular}

\vspace*{.1in}

Statistics Technical Report \\
October 1, 2009

\vspace*{.1in}

\begin{abstract}
Consider the standard linear regression model $\y = \Xmat \betastar +
w$, where $\y \in \real^\numobs$ is an observation vector, $\Xmat \in
\real^{\numobs \times \pdim}$ is a design matrix, $\betastar \in
\real^\pdim$ is the unknown regression vector, and $w \sim
\mathcal{N}(0, \sigma^2 I)$ is additive Gaussian noise. This paper
studies the minimax rates of convergence for estimation of $\betastar$
for $\ell_\rpar$-losses and in the $\ell_2$-prediction loss, assuming
that $\betastar$ belongs to an $\ell_{\qpar}$-ball $\Ballq(\myrad)$
for some $\qpar \in [0,1]$.  We show that under suitable regularity
conditions on the design matrix $\Xmat$, the minimax error in
$\ell_2$-loss and $\ell_2$-prediction loss scales as $\Rq \big(\frac{\log
\pdim}{n}\big)^{1-\frac{\qpar}{2}}$. In addition, we provide lower
bounds on minimax risks in $\ell_{\rpar}$-norms, for all $\rpar \in
[1, +\infty], \rpar \neq \qpar$.  Our proofs of the lower bounds are
information-theoretic in nature, based on Fano's inequality and
results on the metric entropy of the balls $\Ballq(\myrad)$, whereas
our proofs of the upper bounds are direct and constructive, involving
direct analysis of least-squares over $\ell_{\qpar}$-balls.  For the
special case $q = 0$, a comparison with $\ell_2$-risks achieved by
computationally efficient $\ell_1$-relaxations reveals that although
such methods can achieve the minimax rates up to constant factors,
they require slightly stronger assumptions on the design matrix
$\Xmat$ than algorithms involving least-squares over the
$\ell_0$-ball.
\end{abstract}
\end{center}

}
{}
%%%%%%%% END STAT TECHNICAL REPORT FRONT PAGE %%%%%%%%%%%%%%%%

\section{Introduction}

The area of high-dimensional statistical inference concerns the
estimation in the ``large $\pdim$, small $\numobs$'' regime, where
$\pdim$ refers to the ambient dimension of the problem and $\numobs$
refers to the sample size.  Such high-dimensional inference problems
arise in various areas of science, including astrophysics, remote
sensing and geophysics, and computational biology, among others. In
the absence of additional structure, it is frequently impossible to
obtain consistent estimators unless the ratio $\pdim/\numobs$
converges to zero.  However, many applications require solving
inference problems with $\pdim \geq \numobs$, so that consistency is
not possible without imposing additional structure.  Accordingly, an
active line of research in high-dimensional inference is based on
imposing various types of structural conditions, such as sparsity,
manifold structure, or graphical model structure, and then studying
the performance of different estimators.  For instance, in the case of
models with some type of sparsity constraint, a great deal of of work
has studied the behavior of $\ell_1$-based relaxations.

Complementary to the understanding of computationally efficient
procedures are the fundamental or information-theoretic limitations of
statistical inference, applicable to any algorithm regardless of its
computational cost.  There is a rich line of statistical work on such
fundamental limits, an understanding of which can have two types of
consequences.  First, they can reveal gaps between the performance of
an optimal algorithm compared to known computationally efficient
methods.  Second, they can demonstrate regimes in which practical
algorithms achieve the fundamental limits, which means that there is
little point in searching for a more effective algorithm.  As we shall
see, the results in this paper lead to understanding of both types.

\subsection{Problem set-up}

The focus of this paper is a canonical instance of a high-dimensional
inference problem, namely that of linear regression in $\pdim$
dimensions with sparsity constraints on the regression vector
$\betastar \in \real^\pdim$.  In this problem, we observe a pair $(\y,
\Xmat) \in \real^\numobs \times \real^{\numobs \times \pdim}$, where
$\Xmat$ is the design matrix and $\y$ is a vector of response
variables.  These quantities are linked by the standard linear model
\begin{eqnarray}
\label{EqnLinearObs}
\y & = & \Xmat \betastar + \w, 
\end{eqnarray}
where $w \sim N(0, \sigma^2 I_{\numobs \times \numobs})$ is
observation noise.  The goal is to estimate the unknown vector
$\betastar \in \real^\pdim$ of regression coefficients. The sparse
instance of this problem, in which $\betastar$ satisfies some type of
sparsity constraint, has been investigated extensively over the past
decade. Let $X_i$ denote the $i^{th}$ row of $\Xmat$ and $X_j$ denote
the $j^{th}$ column of $\Xmat$. A variety of practical algorithms have
been proposed and studied, many based on $\ell_1$-regularization,
including basis pursuit~\cite{Chen98}, the Lasso~\cite{Tibshirani96},
and the Dantzig selector~\cite{CandesTao06}.  Various authors have
obtained convergence rates for different error metrics, including
$\ell_2$-error~\cite{CandesTao06,BicRitTsy08,Huang06}, prediction
loss~\cite{BicRitTsy08, GreenRitov04}, as well as model selection
consistency~\cite{Huang06,Meinshausen06, Wainwright06,Zhao06}.  In
addition, a range of sparsity assumptions have been analyzed,
including the case of \emph{hard sparsity} meaning that $\betastar$
has exactly $\s \ll \pdim$ non-zero entries, or \emph{soft sparsity}
assumptions, based on imposing a certain decay rate on the ordered
entries of $\betastar$.

\paragraph{Sparsity constraints} These notions of sparsity 
can be defined more precisely in terms of the
$\ell_{\qpar}$-balls\footnote{Strictly speaking, these sets are not
``balls'' when $\qpar < 1$, since they fail to be convex.}  for $\qpar
\in [0,1]$, defined as
\begin{eqnarray}
\Ballq(\myrad) & \defn & \big \{ \beta \in \real^\pdim \, \mid \,
\|\beta\|_{\qpar}^{\qpar} = \sum_{j=1}^\pdim |\beta_j|^\qpar \leq \Rq \big \},
\end{eqnarray}
where in the limiting case $\qpar = 0$, we have the
$\ell_0$-ball
\begin{eqnarray}
\BZ & \defn & \big \{ \beta \in \real^\pdim \, \mid \,
\sum_{j=1}^\pdim \Ind[\beta_j \neq 0] \leq \s \big \},
\end{eqnarray}
corresponding to the set of vectors $\beta$ with at most
$\s$ non-zero elements.

\paragraph{Loss functions} We consider estimators 
$\betahat: \real^\numobs \times \real^{\numobs \times \pdim}
\rightarrow \real^\pdim$ that are measurable functions of the data
$(y, \Xmat)$.  Given any such estimator of the true parameter
$\betastar$, there are many criteria for determining the quality of
the estimate.  In a decision-theoretic framework, one introduces a
loss function such that $\Loss(\betahat, \betastar)$ represents the
loss incurred by estimating $\betahat$ when $\betastar \in \Bq$ is the
true parameter.  The associated risk $\Risk$ is the expected value of
the loss over distributions of ($Y$, $\Xmat$)---namely, the quantity
$\Risk(\betahat, \betastar) = \Exs[\Loss(\betahat, \betastar)]$.
Finally, in the minimax formalism, one seeks to choose an estimator
that minimizes the worst-case risk given by
\begin{equation}
\label{EqnDefnMinmax}
\min_{\betahat} \max_{\betastar \in \Bq} \Risk(\betahat, \betastar).
\end{equation}

Various choices of the loss function are possible, including (a) the
\emph{model selection loss}, which is zero if
$\mbox{supp}(\betahat)=\mbox{supp}(\betastar)$ and one otherwise; (b)
the \emph{$\ell_\rpar$-losses}
\begin{eqnarray}
\label{EqnRparLoss}
\Loss_{\rpar}(\betahat, \betastar) & \defn & \|\betahat -
\betastar\|_\rpar^\rpar \; = \; \sum_{j=1}^\pdim |\betahat_j -
\betastar_j|_\rpar^\rpar,
\end{eqnarray}
and (c) the $\LTwo$-\emph{prediction loss} $\|\Xmat(\betahat -
\betastar)\|_2^2/\numobs$.  In this paper, we study the
$\ell_\rpar$-losses and the $\LTwo$-prediction loss.

\subsection{Our main contributions and related work}

In this paper, we study minimax risks for the high-dimensional linear
model~\eqref{EqnLinearObs}, in which the regression vector $\betastar$
belongs to the ball $\Ballq(\myrad)$ for $0 \leq \qpar \leq 1$. The
core of the paper consists of four main theorems, corresponding to
lower bounds on minimax rate for the cases of $\ell_\rpar$ losses and
the $\LTwo$-prediction loss, and upper bounds for $\ell_2$-norm loss
and the $\LTwo$-prediction loss.  More specifically, in
Theorem~\ref{ThmLower}, we provide lower bounds for
$\ell_{\rpar}$-losses that involve a maximum of two quantities: a term
involving the diameter of the null-space restricted to the
$\ell_{\qpar}$-ball, measuring the degree of non-identifiability of
the model, and a term arising from the $\ell_{\rpar}$-metric entropy
structure for $\ell_{\qpar}$-balls, measuring the massiveness of the
parameter space. Theorem~\ref{ThmUpper} is complementary in nature,
devoted to upper bounds for $\ell_2$-loss. For $\ell_2$-loss, the
upper and lower bounds match up to factors independent of the triple
$(\numobs, \pdim, \myrad)$, and depend only on structural properties
of the design matrix $\Xmat$ (see Theorems~\ref{ThmLower}
and~\ref{ThmUpper}).  Finally, Theorems~\ref{ThmLowerPred}
and~\ref{ThmUpperPred} provide upper and lower bounds for $\LTwo$-prediction
loss. For the $\LTwo$-prediction loss, we provide upper and lower
bounds on minimax risks that are again matching up to factors
independent of $(\numobs, \pdim, \myrad)$, as summarized in
Theorems~\ref{ThmLowerPred} and~\ref{ThmUpperPred}. Structural
properties of the design matrix $\Xmat$ again play a role in minimax
$\LTwo$-prediction risks, but enter in a rather different way than in
the case of $\ell_2$-loss.

For the special case of the Gaussian sequence model (where $\Xmat =
\sqrt{\numobs} I_{n \times n}$), our work is closely related to the
seminal work by of Donoho and Johnstone~\cite{DonJoh94}, who
determined minimax rates for $\ell_\rpar$-losses over
$\ell_\qpar$-balls.  Our work applies to the case of general $\Xmat$,
in which the sample size $\numobs$ need not be equal to the dimension
$\pdim$; however, we re-capture the same scaling as Donoho and
Johnstone~\cite{DonJoh94} when specialized to the case $\Xmat =
\sqrt{\numobs} I_{\numobs \times \numobs}$.  In addition to our analysis of
$\ell_\rpar$-loss, we also determine minimax rates for
$\LTwo$-prediction loss which, as mentioned above, can behave very
differently from the $\ell_2$-loss for general design matrices
$\Xmat$.  During the process of writing up our results, we became
aware of concurrent work by Zhang (see the brief
report~\cite{Zhang09}) that also studies the problem of determining
minimax upper and lower bounds for $\ell_\rpar$-losses with
$\ell_\qpar$-sparsity.  We will be able to make a more thorough
comparison once a more detailed version of their work is publicly
available.

Naturally, our work also has some connections to the vast body of work
on $\ell_1$-based methods for sparse estimation, particularly for the
case of hard sparsity ($q = 0$).  Based on our results, the rates that
are achieved by $\ell_1$-methods, such as the Lasso and the Dantzig
selector, are minimax optimal for $\ell_2$-loss, but require somewhat
stronger conditions on the design matrix than an ``optimal''
algorithm, which is based on searching the $\ell_0$-ball.  We compare
the conditions that we impose in our minimax analysis to various
conditions imposed in the analysis of $\ell_1$-based methods,
including the restricted isometry property of Candes and
Tao~\cite{CandesTao06}, the restricted eigenvalue condition imposed in
Menshausen and Yu~\cite{MeinYu06}, the partial Riesz condition in
Zhang and Huang~\cite{Huang06} and the restricted eigenvalue condition
of Bickel et al.~\cite{BicRitTsy08}. We find that ``optimal'' methods,
which are based on minimizing least-squares directly over the
$\ell_0$-ball, can succeed for design matrices where $\ell_1$-based
methods are not known to work.

The remainder of this paper is organized as follows.  In
Section~\ref{SecMain}, we begin by specifying the assumptions on the
design matrix that enter our analysis, and then state our main
results.  Section~\ref{SecConsequences} is devoted to discussion of
the consequences of our main results, including connections to the
normal sequence model, Gaussian random designs, and related results on
$\ell_1$-based methods.  In Section~\ref{SecProof}, we provide the
proofs of our main results, with more technical aspects deferred to
the appendices.

\section{Main results}
\label{SecMain}

This section is devoted to the statement of our main results, and
discussion of some of their consequences.  We begin by specifying the
conditions on the high-dimensional scaling and the design matrix
$\Xmat$ that enter different parts of our analysis, before giving
precise statements of our main results.

In this paper, our primary interest is the high-dimensional regime in
which $\pdim \gg \numobs$.  For technical reasons, for $\qpar \in
(0,1]$, we require the following condition on the scaling of
$(\numobs, \pdim, \Rq)$:
\begin{equation}
\frac{\pdim}{\Rq \numobs^{\qpar/2}} = \Omega (\pdim^\kappa) \quad
\mbox{for some $\kappa > 0$.}
\end{equation}
In the regime $\pdim \geq \numobs$, this assumption will be satisfied
for all $\qpar \in (0,1]$ as long as $\Rq = o(\pdim^{\frac{1}{2} -
\kappa'})$ for some $\kappa' \in (0, 1/2)$, which is a reasonable
condition on the radius of the $\ell_\qpar$-ball for sparse models.
In the work of Donoho and Johnstone~\cite{DonJoh94} on the normal
sequence model (special case of $\Xmat = I$), discussed at more length
in the sequel, the effect of the scaling of the quantity
$\frac{\pdim}{\Rq \numobs^{\qpar/2}}$ on the rate of convergence also
requires careful treatment.

\subsection{Assumptions on design matrices}

Our first assumption, imposed throughout all of our analysis, is that
the columns \mbox{$\{\Xmat_j, j = 1, \ldots, \pdim \}$} of the design
matrix $\Xmat$ are bounded in $\ell_2$-norm:

\bass[Column normalization]
\label{AssColNorm}
There exists a constant $0 < \colnorm < +\infty$ such that
\begin{eqnarray}
\frac{1}{\sqrt{\numobs}} \max_{j=1, \ldots, \pdim} \| \Xmat_j\|_2 &
\leq & \colnorm.
\end{eqnarray}
\eass
In addition, some of our results involve the set defined by
intersecting the kernel of $\Xmat$ with the $\ell_\qpar$-ball, which
we denote $\Kx \defn \mykernel(\Xmat) \cap \Ballq(\myrad)$.  We define
the \emph{$\Ballq(\myrad)$-kernel diameter} in the $\ell_\rpar$-norm
\begin{eqnarray}
\label{EqnDefnDiamKx}
\diam_\rpar(\Kx) & \defn & \max_{\mytheta \in \Kx} \|\mytheta\|_\rpar
\; = \; \max_{\|\mytheta\|_\qpar^\qpar \leq \myrad, \; \Xmat \mytheta
= 0} \|\mytheta\|_\rpar.
\end{eqnarray}
The significance of this diameter should be apparent: for any
``perturbation'' $\Delta \in \Kx$, it follows immediately from the
linear observation model~\eqref{EqnLinearObs} that no method could
ever distinguish between $\betastar = 0$ and $\betastar = \Delta$.
Consequently, this $\Ballq(\myrad)$-kernel diameter is a measure of
the \emph{lack of identifiability} of the linear
model~\eqref{EqnLinearObs} over $\Ballq(\myrad)$.

Our second assumption, which is required only for achievable results
for $\ell_2$-error and lower bounds for $\LTwo$-prediction error,
imposes a lower bound on $\|\Xmat \mytheta\|_2/\sqrt{\numobs}$ in
terms of $\|\mytheta\|_2$ and a residual term:

\bass[Lower bound on restricted curvature]
\label{AssLower}
There exists a constant \mbox{$\conlower > 0$} and a function
$\SpecFunLower$ such that
\begin{eqnarray}
\label{EqnAssLower}
\frac{1}{\sqrt{\numobs}} \| \Xmat \mytheta \|_2 & \geq &
\conlower \, \| \mytheta \|_2 - \SpecFunLower
\mbox{\;\;for all \;\;$\mytheta \in \Ballq(2 \myrad)$.}
\end{eqnarray}
\eass
\vspace*{.1in}

\noindent {\bf{Remarks:}} Conditions on the scaling for $\SpecFunLower$ are provided in Theorems~\ref{ThmUpper} and \ref{ThmLowerPred}. It is useful to recognize that the lower bound~\eqref{EqnAssLower} is closely related to the diameter condition~\eqref{EqnDefnDiamKx}; in particular,
Assumption~\ref{AssLower} induces an upper bound on the
$\Ballq(\myrad)$-kernel diameter in $\ell_2$-norm, and hence the
identifiability of the model:
\blems
\label{LemIdent}
If Assumption~\ref{AssLower} holds for any $\qpar \in (0,1]$, then the
$\Ballq(\myrad)$-kernel diameter in $\ell_2$-norm is upper bounded as
\begin{eqnarray*}
\diam_2(\Kx) & \leq & \frac{\SpecFunLower}{\conlower}.
\end{eqnarray*}
\elems
\spro
We prove the contrapositive statement. Note that if $\diam_2(\Kx) >
\frac{\SpecFunLower}{\conlower}$, then there must exist some $\mytheta
\in \Ballq(\myrad)$ with $\Xmat \mytheta = 0$ and $\|\mytheta\|_2 \; >
\: \frac{\SpecFunLower}{\conlower}$. We then conclude that
\begin{eqnarray*}
0 \; = \; \frac{1}{\sqrt{\numobs}} \|\Xmat \mytheta\|_2 & < &
\conlower \|\mytheta\|_2 - \SpecFunLower,
\end{eqnarray*}
which implies there cannot exist any $\conlower$ for which the
lower bound~\eqref{EqnAssLower} holds.  
\fpro 

In Section~\ref{SecCompareEllone}, we discuss further connections
between our assumptions, and the conditions imposed in analysis of the
Lasso and other $\ell_1$-based
methods~\cite{CandesTao06,Meinshausen06,BicRitTsy08}.  In the case
$\qpar = 0$, we find that Assumption~\ref{AssLower} is weaker than any
condition under which an $\ell_1$-based method is known to succeed.
Finally, in Section~\ref{SecRandDesign}, we prove that versions of
both Assumptions~\ref{AssColNorm} and~\ref{AssLower} hold with high
probability for various classes of non-i.i.d. Gaussian random design
matrices (see Proposition~\ref{PropRandDesign}).

\subsection{Risks in $\ell_\rpar$-norm}
\label{SecRparNorm}

Having described our assumptions on the design matrix, we now turn to
the main results that provide upper and lower bounds on minimax
risks. In all of the statements to follow, we use the quantities
$\AnnoyLower$, $\AnnoyUpper$, $\AnnoyPred$ etc. to denote numerical
constants, independent of $\numobs, \pdim$, $\myrad$, $\sigma^2$ and
the design matrix $\Xmat$.  We begin with lower bounds on the
$\ell_\rpar$-risk.
\btheos[Lower bounds on $\ell_\rpar$-risk]
\label{ThmLower}
Consider the linear model~\eqref{EqnLinearObs} for a fixed design
matrix $\Xmat \in \real^{\numobs \times \pdim}$. 
\begin{enumerate}
\item[(a)] {\bf{Conditions for $\qpar \in (0,1]$:}} Suppose that
$\Xmat$ is column-normalized (Assumption~\ref{AssColNorm} with
$\colnorm < \infty$).  For any $\rpar \in [1, \infty)$, the minimax
$\ell_\rpar$-risk over the $\ell_\qpar$ ball is lower bounded as
\begin{eqnarray}
\label{EqnLowerBound} 
\min_{\betahat} \max_{\betastar \in \Bq} \mathbb{E}\|\betahat -
\betastar\|_{\rpar}^\rpar & \geq & \AnnoyLower \; \max \Biggr \{
\diam_{\rpar}^{\rpar}(\Kx), \: \myrad \,
\biggr [\frac{\sigma^2}{\colnormsq} \, \frac{\log
\pdim}{\numobs} \biggr]^{\frac{\rpar-\qpar}{2}} \Biggr \}. 
\end{eqnarray}
\item[(b)] {\bf{Conditions for $\qpar = 0$:}} 
Suppose that
$\frac{\|\Xmat \mytheta\|_2}{\sqrt{\numobs} \|\mytheta\|_2} \leq
\conupper$ for all $\mytheta \in \Ball_0(2 \s)$.  Then for any
$\rpar \in [1, \infty)$, the minimax $\ell_\rpar$-risk over the
$\ell_0$-ball with radius $\s = R_0$ is lower bounded as
\begin{eqnarray}  
\label{EqnLowerBoundZero}
\min_{\betahat} \max_{\betastar \in \BZ} \mathbb{E}\|\betahat -
\betastar\|_{\rpar}^{\rpar} & \geq & \AnnoyLowerZero \; \max \Biggr \{
\diam_{\rpar}^{\rpar}(\KZx), \; \spindex^{\frac{\rpar}{2}} \big[
\frac{\sigma^2}{\conuppersq} \frac{\log(\pdim/\spindex)}{n}
\big]^{\frac{\rpar}{2}} \Biggr \}. 
\end{eqnarray}
\end{enumerate}
\etheos

Note that both lower bounds consist of two terms.  The first term is
simply the diameter of the set $\Kx = \mykernel(\Xmat) \cap
\Ballq(\myrad)$, which reflects the extent which the linear
model~\eqref{EqnLinearObs} is unidentifiable.  Clearly, one cannot
estimate $\betastar$ any more accurately than the diameter of this
set.  In both lower bounds, the ratios $\sigma^2/\colnormsq$
(or $\sigma^2/\conuppersq$) correspond to the inverse of the
signal-to-noise ratio, comparing the noise variance $\sigma^2$ to the
magnitude of the design matrix measured by $\conupper$. As the proof will clarify, the term $[\log \pdim]^{\frac{\rpar - \qpar}{2}}$ in the lower bound~\eqref{EqnLowerBound}, and similarly the term $\log(\frac{\pdim} {\s})$ in the bound~\eqref{EqnLowerBoundZero}, are
reflections of the complexity of the $\ell_\qpar$-ball, as measured by
its metric entropy.  For many classes of design matrices, the second
term is of larger order than the diameter term, and hence determines
the rate.  (In particular, see Section~\ref{SecRandDesign} for an
in-depth discussion of the case of random Gaussian designs.)

\vspace*{.2in}

We now state upper bounds on the $\ell_2$-norm minimax risk over
$\ell_\qpar$ balls.  For these results, we require both the column
normalization condition (Assumption~\ref{AssColNorm}) and the
curvature condition (Assumption~\ref{AssLower}).

\btheos [Upper bounds on $\ell_2$-risk]
\label{ThmUpper}
Consider the model~\eqref{EqnLinearObs} with a fixed design matrix
\mbox{$\Xmat \in \real^{\numobs \times \pdim}$} that is
column-normalized (Assumption~\ref{AssColNorm} with $\colnorm <
\infty$).
\begin{enumerate}
\item[(a)] {\bf{Conditions for $\qpar \in (0,1]$:}} If $\Xmat$
satisfies Assumption~\ref{AssLower} with \mbox{$\SpecFunLower =
o(\Rq^{1/2}(\frac{\log \pdim}{n})^{1/2-\qpar/4})$} and $\conlower >
0$, then there exist constants $c_1$ and $c_2$ such that the minimax
$\ell_2$-risk is upper bounded as
\begin{eqnarray}
\label{EqnUpperBound}
 \min_{\betahat} \max_{\betastar \in \Bq} \|\betahat - \betastar\|_2^2
& \leq & 24 \myrad \Big[ \frac{\colnormsq}{\conlowersq}
\frac{\sigma^2}{\conlowersq} \; \frac{\log \pdim}{\numobs} \Big]^{1 -
\qpar/2},
\end{eqnarray}
with probability greater than $1 - c_1\exp{(-c_2 \numobs)}$.
\item[(b)] {\bf{Conditions for $\qpar = 0$:}} If $\Xmat$ satisfies
Assumption~\ref{AssLower} with $\SpecFunLowerZ = 0$ and $\conlower >
0$, then there exists constants $c_1$ and $c_2$ such that the minimax
$\ell_2$-risk is upper bounded as
\begin{eqnarray}
\label{EqnUpperBoundZero}
 \min_{\betahat} \max_{\betastar \in \BZ} \|\betahat -
\betastar\|_2^2 & \leq & 6
\frac{\colnormsq}{\conlowersq} \;
\frac{\sigma^2}{\conlowersq} \; \frac{\spindex \, \log
\pdim}{\numobs}, \qquad
\end{eqnarray}
 with probability greater than $1 - c_1\exp{(-c_2 n)}$.  If, in
addition, the design matrix satisfies $\frac{\|\Xmat
\mytheta\|_2}{\sqrt{\numobs} \|\mytheta\|_2} \leq \conupper$ for all
$\mytheta \in \Ball_0(2 \s)$, then the minimax $\ell_2$-risk is upper
bounded as
\begin{eqnarray}
\label{EqnUpperBoundZeroSharp}
 \min_{\betahat} \max_{\betastar \in \BZ} \|\betahat - \betastar\|_2^2
& \leq & 144 \: \frac{\conuppersq}{\conlowersq}
\frac{\sigma^2}{\conlowersq} \, \frac{\spindex
\log(\pdim/\spindex)}{\numobs},
\end{eqnarray}
with probability greater than $1 - c_1\exp{(-c_2 \s \log(\pdim -\s)
  )}$.
\end{enumerate}
\etheos 

In the case of $\ell_2$-risk and design matrices $\Xmat$ that satisfy
the assumptions of both Theorems~\ref{ThmLower} and~\ref{ThmUpper},
then these results identify the minimax risk up to constant factors.
In particular, for $\qpar \in (0,1]$, the minimax $\ell_2$-risk scales
as
\begin{equation}
 \min_{\betahat} \max_{\betastar \in \Ballq(\myrad)} \Exs \|\betahat -
\betastar\|_2^2 \; = \; \Theta \Big(\myrad \Big[\frac{\sigma^2 \log
\pdim}{\numobs} \Big]^{1 - \qpar/2} \Big),
\end{equation}
whereas for $\qpar = 0$, the minimax $\ell_2$-risk scales as
\begin{equation}
 \min_{\betahat} \max_{\betastar \in \BZ} \mathbb{E}\|\betahat -
\betastar\|_2^2 \; = \; \Theta \Big(\frac{\sigma^2 \, \spindex
\log(\pdim/\spindex)}{\numobs} \Big).
\end{equation}
Note that the bounds with high probability can be converted to bound
in expectation by a standard integration over the tail probability.

\subsection{Risks in prediction norm}

In this section, we investigate minimax risks in terms of the
$\LTwo$-prediction loss $\|\Xmat(\betahat -
\betastar)\|_2^2/\numobs$, and provide both lower and upper
bounds on it.
\btheos[Lower bounds on prediction risk]
\label{ThmLowerPred}
Consider the model~\eqref{EqnLinearObs} with a fixed design matrix
$\Xmat \in \real^{\numobs \times \pdim}$ that is column-normalized
(Assumption~\ref{AssColNorm} with $\colnorm < \infty$).
\begin{enumerate}
\item[(a)] {\bf{Conditions for $\qpar \in (0,1]$:}} If the design
matrix $\Xmat$ satisfies Assumption~\ref{AssLower} with $\conlower >
0$ and \mbox{$\SpecFunLower = o(\Rq^{1/2}(\frac{\log
\pdim}{n})^{1/2-\qpar/4})$,} then the minimax prediction risk is lower
bounded as
\begin{eqnarray}
\label{EqnLowerPred}
\min_{\betahat} \max_{\beta \in \Bq} \Exs \frac{\|\Xmat \, (\betahat -
\beta) \|_2^2}{\numobs} & \geq & \Annoy'_{2, \qpar} \; \Rq \;
\conlowersq \; \Big[\frac{\sigma^2}{\colnormsq} \; \frac{\log
\pdim}{\numobs} \Big]^{1 - \qpar/2}.
\end{eqnarray}
\item[(b)] {\bf{Conditions for $\qpar = 0$:}} Suppose that $\Xmat$
satisfies Assumption~\ref{AssLower} with $\conlower > 0$ and
\mbox{$\SpecFunLowerZ = 0$,} and moreover that $\frac{\|\Xmat
\mytheta\|_2}{\sqrt{\numobs} \|\mytheta\|_2} \leq \conupper$ for all
$\mytheta \in \Ball_0(2 \s)$.  Then the minimax prediction risk is
lower bounded as
\begin{eqnarray}
\label{EqnLowerPredZero}
\min_{\betahat} \max_{\beta \in \BZ} \Exs \frac{\|\Xmat (\betahat -
\beta)\|_2^2}{\numobs} & \geq & \Annoy'_{0, \qpar} \; \conlowersq \;
\frac{\sigma^2}{\conuppersq} \; \frac{\spindex \log (\pdim/\spindex)}{n}.
\end{eqnarray}
\end{enumerate}
\etheos

In the other direction, we have the following result:
\btheos[Upper bounds on prediction risk]
\label{ThmUpperPred}
Consider the model~\eqref{EqnLinearObs} with a fixed design matrix
$\Xmat \in \real^{\numobs \times \pdim}$.  
\begin{enumerate}
\item[(a)] {\bf{Conditions for $\qpar \in (0,1]$:}} If $\Xmat$
satisfies the column normalization condition, then for some constant
$\Annoy_{2, \qpar}$, there exist $c_1$ and $c_2$ such that the minimax
prediction risk is upper bounded as
\begin{eqnarray}
\label{EqnUpperPred} 
\min_{\betahat} \max_{\betastar \in \Bq} \frac{1}{\numobs}\|\Xmat(\betahat - \betastar)\|_2^2 &
\leq & \Annoy_{2, \qpar} \; \colnormsq \; \myrad \;
\Big[\frac{\sigma^2}{\colnormsq} \; \frac{\log \pdim}{\numobs}
\Big]^{1-\frac{\qpar}{2}},
\end{eqnarray}
with probability greater than $1 - c_1 \exp{(-c_2 \Rq (\log
\pdim)^{1-\qpar/2}n^{\qpar/2})}$.
\item[(b)] {\bf{Conditions for $\qpar = 0$:}} For any $\Xmat$, with
probability greater than $1 - \exp{(-10\s \log (\pdim/\s))}$ the
minimax prediction risk is upper bounded as
\begin{eqnarray}
\label{EqnUpperPredZero}
 \min_{\betahat} \max_{\betastar \in \BZ} \frac{1}{\numobs} \| \Xmat(\betahat - \betastar) \|_2^2 &
 \leq & 81 \frac{\sigma^2 \; \spindex \log(\pdim/\spindex)}{\numobs}.
\end{eqnarray}
\end{enumerate}
\etheos 

\subsection{Some intuition}

In order to provide the reader with some intuition, let us make some
comments about the scalings that appear in our results.

First, as a basic check of our results, it can be verified that
Lemma~\ref{LemIdent} ensures that the lower bounds on minimax rates
stated in Theorem~\ref{ThmLower} for $\rpar = 2$ are always less than
or equal to the achievable rates stated in Theorem~\ref{ThmUpper}.  In
particular, since $\SpecFunLower = o(\Rq^{1/2}(\frac{\log
\pdim}{n})^{1/2-\qpar/4})$ for $\qpar \in (0,1]$, Lemma~\ref{LemIdent}
implies that $\diam_2^2(\Kx) = o(\Rq(\frac{\log
\pdim}{n})^{1-\qpar/2})$, meaning that the achievable rates are always
at least as large as the lower bounds in the case $q \in (0,1]$.  In
the case of hard sparsity ($\qpar = 0$), the upper and lower bounds
are clearly consistent since $\SpecFunLowerZ = 0$ implies the diameter
of $\KZx$ is $0$.

Second, for the case $\qpar=0$, there is a concrete interpretation of
the rate $\frac{\s \log (\pdim/\s)}{n}$, which appears in
Theorems~\ref{ThmLower}(b), ~\ref{ThmUpper}(b), ~\ref{ThmLowerPred}(b)
and~\ref{ThmUpperPred}(b)).  Note that there are ${\pdim \choose \s}$
subsets of size $\s$ within $\{ 1, 2, \ldots, \pdim \}$, and by
standard bounds on binomial coefficients~\cite{Cover}, we have $\log
{\pdim \choose \s} = \Theta(\s \log(\pdim/\s))$.  Consquently, the
rate $\frac{\s \log (\pdim/\s)}{n}$ corresponds to the log number of
models divided by the sample size $\numobs$.  Note that unless
$\s/\pdim = \Theta(1)$, this rate is equivalent (up to constant
factors) to $\frac{ \s \log \pdim}{\numobs}$.

Third, for $\qpar \in (0, 1]$, the interpretation of the rate $\Rq
\big(\frac{\log \pdim}{n}\big)^{1-\qpar/2}$, appearing in parts (a) of
Theorems~\ref{ThmLower} through~\ref{ThmUpperPred}, is less
immediately obvious but can can understood as follows.  Suppose that
we choose a subset of size $\s_{\qpar}$ of coefficients to estimate,
and ignore the remaining $\pdim - s_{\qpar}$ coefficients.  For
instance, if we were to choose the top $\s_\qpar$ coefficients of
$\beta^*$ in absolute value, then the fast decay imposed by the
$\ell_\qpar$-ball condition on $\beta^*$ would mean that the remaining
$\pdim - \s_\qpar$ coefficients would have relatively little impact.
With this intuition, the rate for $\qpar > 0$ can be interpreted as
the rate that would be achieved by choosing $\s_{\qpar} = \Rq
\big(\frac{\log \pdim}{n}\big)^{-\qpar/2}$, and then acting as if the
problem were an instance of a hard-sparse problem ($\qpar = 0$) with
$s = \s_\qpar$.  For such a problem, we would expect to achieve the
rate $\frac{s_\qpar \log \pdim}{\numobs}$, which is exactly equal to
$\Rq \big(\frac{\log \pdim}{n}\big)^{1-\qpar/2}$.  Of course, we have
only made a very heuristic argument here; this truncation idea is made
more precise in Lemma~\ref{LemElloneElltwo} to appear in the sequel.

Fourth, we note that the minimax rates for $\LTwo$-prediction error
and $\ell_2$-norm error are essentially the same except that the
design matrix structure enters minimax risks in \emph{very different
ways}.  In particular, note that proving lower bounds on prediction
risk requires imposing relatively strong conditions on the design
$\Xmat$---namely, Assumptions~\ref{AssColNorm} and~\ref{AssLower} as
stated in Theorem~\ref{ThmLowerPred}. In contrast, obtaining upper
bounds on prediction risk requires very mild conditions. At the most
extreme, the upper bound for $\qpar = 0$ in Theorem~\ref{ThmLowerPred}
requires no assumptions on $\Xmat$ while for $\qpar > 0$ only the
column normalization condition is required. All of these statements
are reversed for $\ell_2$-risks, where lower bounds can be proved with
only Assumption~\ref{AssColNorm} on $\Xmat$ (see
Theorem~\ref{ThmLower}), whereas upper bounds require both
Assumptions~\ref{AssColNorm} and~\ref{AssLower}.

Lastly, in order to appreciate the difference between the conditions
for $\LTwo$-prediction error and $\ell_{2}$ error, it is useful to
consider a toy but illuminating example.  Consider the linear
regression problem defined by a design matrix $\Xmat =
\begin{bmatrix} \Xmat_1 & \Xmat_2 & \cdots & \Xmat_\pdim
\end{bmatrix}$ with \emph{identical columns}---that is, $X_j =
\widetilde{\Xmat}_1$ for all $j = 1, \ldots, \pdim$.  We assume that
vector $\widetilde{\Xmat}_1 \in \real^\pdim$ is suitably scaled so
that the column-normalization condition (Assumption~\ref{AssColNorm})
is satisfied.  For this particular choice of design matrix, the linear
observation model~\eqref{EqnLinearObs} reduces to $\y =
(\sum_{j=1}^{\pdim}{\beta^*_j}) \widetilde{\Xmat}_1 + \w$.  For the
case of hard sparsity ($\qpar=0$), an elementary argument shows that
the minimax risk in $\LTwo$-prediction error scales as
$\Theta(\frac{1}{\numobs})$.  This scaling implies that the upper
bound~\eqref{EqnUpperPredZero} from Theorem~\ref{ThmUpperPred} holds
(but is not tight).\footnote{Note that the lower
bound~\eqref{EqnLowerPredZero} on the $\LTwo$-prediction error from
Theorem~\ref{ThmLowerPred} does not apply to this model, since this
degenerate design matrix with identical columns does not satisfy any
version of Assumption~\ref{AssLower}.}  Consequently, this highly
degenerate design matrix yields a very easy problem for
$\LTwo$-prediction, since the $1/\numobs$ rate is essentially
parametric.  In sharp contrast, for the case of $\ell_{2}$-norm error
(still with hard sparsity $\qpar = 0$), the model becomes
unidentifiable.  To see the lack of identifiability, let $e_i \in
\real^\pdim$ denote the unit-vector with $1$ in position $i$, and
consider the two regression vectors $\betastar = c \,e_1$ and
$\widetilde{\beta} = c \, e_2$, for some constant $c \in \real$.  Both
choices yield the same observation vector $\y$, and since the choice
of $c$ is arbitrary, the minimax $\ell_2$-error is infinite.  In this
case, the lower bound~\eqref{EqnLowerBoundZero} on $\ell_2$-error from
Theorem~\ref{ThmLower} holds (and is tight, since the kernel diameter
is infinite).  In contrast, the upper bound~\eqref{EqnUpperBoundZero}
on $\ell_2$-error from Theorem~\ref{ThmUpper} does not apply,
because Assumption~\ref{AssLower} is violated due to the extreme
degeneracy of the design matrix.

%%%%%%%%%%%%%%%%%%%%%%%%%%%%%%%%%%%%%%%%%%%%%%%%%%%%%%%%%%%%%%%%%%%%%%%%%%%
%%%%%%%%%%%%%%%%%%%%%%%%%%%%%%%%%%%%%%%%%%%%%%%%%%%%%%%%%%%%%%%%%%%%%%%%%%%%

\section{Some consequences}
\label{SecConsequences}

In this section, we discuss some consequences of our results.  We
begin by considering the classical Gaussian sequence model, which
corresponds to a special case of our linear regression model, and
making explicit comparisons to the results of Donoho and
Johnstone~\cite{DonJoh94} on minimax risks over $\ell_\qpar$-balls.

\subsection{Connections with the normal sequence model}

The normal (or Gaussian) sequence model is defined by the observation
sequence
\begin{eqnarray}
\label{EqnGaussSeqModel} 
y_i & = & \thetastar_i + \gnoise_i, \qquad \mbox{for $i=1, \ldots,
\numobs$,}
\end{eqnarray}
where $\thetastar \in \ThetaSet \subseteq \real^\numobs$ is a fixed
but unknown vector, and the noise variables $\gnoise_i \sim
\mathcal{N}(0,\frac{\tau^2}{\numobs})$ are i.i.d. normal variates.
Many non-parametric estimation problems, including regression and
density estimation, are asymptotically equivalent to an instance of
the Gaussian sequence model~\cite{Pinsk80,Nuss96,BrownLo96}, where the
set $\ThetaSet$ depends on the underlying ``smoothness'' conditions
imposed on the functions.  For instance, for functions that have an
$m^{th}$ derivative that is square-differentiable (a particular kind
of Sobolev space), the set $\ThetaSet$ corresponds to an ellipsoid; on
the other hand, for certain choices of Besov spaces, it corresponds to
an $\ell_\qpar$-ball.

In the case $\ThetaSet = \Ball_\qpar(\myrad)$, our linear regression
model~\eqref{EqnLinearObs} includes the normal sequence
model~\eqref{EqnGaussSeqModel} as a special case.  In particular, it
corresponds to setting $\pdim = \numobs$, the design matrix $\Xmat =
I_{n \times n}$, and noise variance $\sigma^2 =
\frac{\tau^2}{\numobs}$.  For this particular model, seminal work by
Donoho and Johnstone~\cite{DonJoh94} derived sharp asymptotic results
on the minimax error for general $\ell_\rpar$-norms over $\ell_\qpar$
balls.  Here we show that a corollary of our main theorems yields the
same scaling in the case $\rpar = 2$ and $\qpar \in [0,1]$.

\begin{cors} 
\label{CorNorSeqModel}
Consider the normal sequence model~\eqref{EqnGaussSeqModel} with
$\ThetaSet = \Ball_\qpar(\myrad)$ for some $\qpar \in (0,1]$.  Then
there are constants $\Annoy'_\qpar \leq \Annoy_\qpar$ depending only
on $\qpar$ such that
\begin{equation}
\Annoy'_\qpar (\frac{2 \tau^2 \log
\numobs}{\numobs})^{1-\frac{\qpar}{2}} \; \leq \; \min_{\betahat}
\max_{\betastar \in \Bq} \Exs \|\betahat - \betastar\|_2^2 \; \leq \;
\Annoy_\qpar (\frac{2 \tau^2 \log
\numobs}{\numobs})^{1-\frac{\qpar}{2}}.
\end{equation}
\end{cors}
%%%%%%%%%%%%%%%%%%%%%%%%%%%%%%%%%%%%%%%%%%%%%%%%%%%%%%%%%%%%%%%%%%%%%%%
These bounds follow from our main theorems, via the substitutions
$\numobs=\pdim$, $\sigma^2 = \frac{\tau^2}{\numobs}$, and
\mbox{$\conupper = \conlower = 1$.}  To be clear, Donoho and
Johnstone~\cite{DonJoh94} provide a far more careful analysis that
yields sharper control of the constants than we have provided here.

%Donoho and
%Johnstone~\cite{DonJoh94} provide results only for $\qpar>0$. However
%our results clearly imply that the same bounds hold for $\qpar=0$.

\subsection{Random Gaussian Design}
\label{SecRandDesign}

Another special case of particular interest is that of random Gaussian
design matrices.  A widely studied instance is the standard Gaussian
ensemble, in which the entries of $\Xmat \in \real^{\numobs \times
\pdim}$ are i.i.d. $N(0,1)$ variates. A variety of results are known
for the singular values of random matrices $\Xmat$ drawn from this
ensemble (e.g.,~\cite{Bai,Silverstein,DavSza01}); moreover, some past
work~\cite{Donoho06a,CandesTao06} has studied the behavior of
different $\ell_1$-based methods for the standard Gaussian ensemble,
in which entries $X_{ij}$ are i.i.d. $N(0, 1)$.  In modeling terms,
requiring that all entries of the design matrix $\Xmat$ are i.i.d. is
an overly restrictive assumption, and not likely to be met in
applications where the design matrix cannot be chosen.  Accordingly,
let us consider the more general class of Gaussian random design
matrices $\Xmat \in \real^{\numobs \times \pdim}$, in which the rows
are independent, but there can be arbitrary correlations between the
columns of $\Xmat$.  To simplify notation, we define the shorthand
$\rhomax(\CovMat) \defn \max_{j=1, \ldots, \pdim} \CovMat_{jj}$,
corresponding to the maximal variance of any element of $\Xmat$, and
use $\CovMat^{1/2}$ to denote the symmetric square root of the
covariance matrix.

In this model, each column $\Xmat_j$, $j = 1, \ldots, \pdim$ has
i.i.d. elements.  Consequently, it is an immediate consequence of
standard concentration results for $\chi^2_\numobs$ variates (see
Appendix~\ref{AppChiTail}) that
\begin{align}
\max_{j =1, \ldots, \pdim} \frac{\|\Xmat_j\|_2}{\sqrt{\numobs}} & \leq
\rhomax(\CovMat) \, \big(1 + \sqrt{\frac{32 \log \pdim}{\numobs}}).
\end{align}
Therefore, Assumption~\ref{AssColNorm} holds as long as $\numobs =
\Omega(\log \pdim)$ and $\rhomax(\CovMat)$ is bounded.

Showing that a version of Assumption~\ref{AssLower} holds with high
probability requires more work.  We summarize our findings in the
following result:
\bprops
\label{PropRandDesign}
Consider a random design matrix $\Xmat \in \real^{\numobs \times
\pdim}$ formed by drawing each row $X_i \in \real^\pdim$ i.i.d. from an
$N(0, \CovMat)$ distribution.  Then for some numerical constants
$\plaincon_k \in (0, \infty)$, $k = 1,2$, we have
\begin{eqnarray}
\label{EqnLowerFinal}
\frac{\|X v\|_2}{\sqrt{\numobs}} & \geq & \frac{1}{2} \|\CovMatSqrt
v\|_2 - 6 \, \big(\frac{\rhomax(\CovMat) \log \pdim}{\numobs}
\big)^{1/2} \; \|v\|_1 \quad \mbox{for all $v \in \real^\pdim$}
\end{eqnarray}
with probability $1- \plaincon_1 \, \exp(-\plaincon_2
\numobs)$.
\eprops

\noindent {\bf{Remarks:}} Past work by by Amini and
Wainwright~\cite{AmiWai08} in the analysis of sparse PCA has
established an upper bound analogous to the lower
bound~\eqref{EqnLowerFinal} for the special case $\CovMat = I_{\pdim
\times \pdim}$.  We provide a proof of this matching upper bound for
general $\CovMat$ as part of the proof of
Proposition~\ref{PropRandDesign} in Appendix~\ref{AppPropRandDesign}.
The argument is based on Slepian's lemma~\cite{DavSza01} and its
extension due to Gordon~\cite{Gordon85}, combined with concentration
of Gaussian measure results~\cite{Ledoux01}.  Note that we have made
no effort to obtain sharp leading constants (i.e., the factors $1/2$
and $6$ can easily be improved), but the basic
result~\eqref{EqnLowerFinal} suffices for our purposes.

 Let us now discuss the implications of this result for
Assumption~\ref{AssLower}.  First, in the case $\qpar = 0$, the
bound~\eqref{EqnUpperBoundZero} in Theorem~\ref{ThmUpper} requires
that Assumption~\ref{AssLower} holds with $\SpecFunLowerZ = 0$ for all
$\mytheta \in \BallZ(2 \s)$. To see the connection with
Proposition~\ref{PropRandDesign}, note that if $\mytheta \in \BallZ(2
\s)$, then we have $\|\mytheta\|_1 \leq \sqrt{2 \s} \|\mytheta\|_2$,
and hence
\begin{align*}
\frac{\|X v\|_2}{\sqrt{\numobs}} & \geq \Big \{ \frac{\|\CovMatSqrt v
\|_2}{ 2 \|v\|_2} - 6 \sqrt{2} \big(\frac{\rhomax(\CovMat) \spindex \log
\pdim}{\numobs} \big)^{1/2} \Big \} \|v\|_2.
\end{align*}
Therefore, as long as $\rhomax(\CovMat) < \infty$, $\min_{v \in
\BallZ(2 \s)} \frac{\|\CovMatSqrt v \|_2}{\|v\|_2} > 0$ and $\frac{\s
\log \pdim}{n} = o(1)$, the condition needed for the
bound~\eqref{EqnUpperBoundZero} will be met.

Second, in the case $\qpar \in (0,1]$, Theorem~\ref{ThmUpper}(a)
requires that Assumption~\ref{AssLower} hold with the residual term
\mbox{$\SpecFunLower = o(\Rq^{1/2} \frac{\log
\pdim}{\numobs})^{1/2-\qpar/4}$}.  We claim that
Proposition~\ref{PropRandDesign} guarantees this condition, as long as
$\rhomax(\CovMat) < \infty$ and the minimum eigenvalue of $\CovMat$ is
bounded away from zero.  In order to verify this claim, we require the
following result:
\blems
\label{LemElloneElltwo}
For any vector $\mytheta \in \Ballq(2 \myrad)$ and any positive number
$\tau > 0$, we have
\begin{eqnarray}
\label{EqnElloneElltwo}
\|\mytheta\|_1 & \leq & \sqrt{2 \myrad} \tau^{-\qpar/2} \|\mytheta\|_2
+ 2 \myrad \tau^{1-\qpar}.
\end{eqnarray}
\elems
\noindent Although this type of result is standard
(e.g,~\cite{DonJoh94}), we provide a proof in
Appendix~\ref{AppLemElloneElltwo} for completeness.  In order to
exploit Lemma~\ref{LemElloneElltwo}, let us set $\tau =
\sqrt{\frac{\log \pdim}{\numobs}}$.  With this choice, we can
substitute the resulting bound~\eqref{EqnElloneElltwo} into the lower
bound~\eqref{EqnLowerFinal}, thereby obtaining that
\begin{align*}
\frac{\|X v\|_2}{\sqrt{\numobs}} & \geq \Big \{ \frac{\|\CovMatSqrt v
\|_2}{ 2 \|v\|_2} - 6 \sqrt{2 \, \rhomax(\CovMat)} \; \sqrt{\Rq} \,
\big(\frac{\log \pdim}{\numobs} \big)^{1/2 - \qpar/4} \big \} \|v\|_2
- 2\Rq \rhomax(\CovMat)^{1/2} \big( \frac{\log
\pdim}{\numobs}\big)^{1-\qpar/2}.
\end{align*}
Recalling that the condition $\sqrt{\Rq} \big( \frac{\log \pdim}{n}
\big)^{1/2-\qpar/4} = o(1)$ is required for consistency, we see that
Assumption~\ref{AssLower} holds as long as $\rhomax(\CovMat) <
+\infty$ and the minimum eigenvalue of $\CovMat$ is bounded away from
zero.

Lastly, it is also worth noting that we can also obtain the following
stronger result for the case $\qpar = 0$, in the case that $\min_{v
\in \BallZ(2 \s)} \frac{\|\CovMatSqrt v \|_2}{\|v\|_2} > 0$ and
$\max_{v \in \BallZ(2 \s)} \frac{\|\CovMatSqrt v \|_2}{\|v\|_2} <
\infty$.  If the sparse eigenspectrum is bounded in this way, then as
long as $\numobs > \plaincon_3 \, \s \log(\pdim/\s)$, we have
\begin{equation}
\label{EqnLowerFinalZero}
3 \|\CovMatSqrt v\|_2 \; \geq \; \frac{\|X v\|_2}{\sqrt{\numobs}} \;
\geq \; \frac{1}{2} \|\CovMatSqrt v\|_2 \quad \mbox{ for all $v \in
\Ball_0(2 \s)$}
\end{equation}
with probability greater than $1 - \plaincon_1 \exp(-\plaincon_2
\numobs)$.  This fact follows by applying the union bound over all
${\pdim \choose 2 \spindex}$ subsets of size $2 \s$, combined with
standard concentration results for random matrices (e.g., see Davidson
and Szarek~\cite{DavSza01} for $\CovMat = I$, and
Wainwright~\cite{Wainwright06} for the straightforward extensions to
non-identity covariances).

\subsection{Comparison to $\ell_1$-based methods}

\label{SecCompareEllone}

In addition, it is interesting to compare our minimax rates of
convergence for $\ell_2$-error with known results for $\ell_1$-based
methods, including the Lasso~\cite{Tibshirani96} and the closely
related Dantzig method~\cite{CandesTao06}. Here we discuss only the
case $\qpar=0$ since we are currently unaware of any $\ell_2$-error
bound for $\ell_1$-based methods for $\qpar \in (0,1]$.  For the
Lasso, past work~\cite{Huang06,MeinYu06} has shown that its
$\ell_2$-error is upper bounded by $\frac{\s \log \pdim}{\numobs}$
under sparse eigenvalue conditions. Similarly, Candes and
Tao~\cite{CandesTao06} show the same scaling for the Dantzig selector,
when applied to matrices that satisfy the more restrictive RIP
conditions.  More recent work by Bickel et. al~\cite{BicRitTsy08}
provides a simultaneous analysis of the Lasso and Dantzig selector
under a common set of assumptions that are weaker than both the RIP
condition and sparse eigenvalue conditions. Together with our results
(in particular, Theorem~\ref{ThmLower}(b)), this body of work shows
that under appropriate conditions on the design $\Xmat$, the rates
achieved by $\ell_1$-methods in the case of hard sparsity ($\qpar =
0$) are minimax-optimal.

Given that the rates are optimal, it is appropriate to compare the
conditions needed by an ``optimal'' algorithm, such as that analyzed
in Theorem~\ref{ThmUpper}, to those used in the analysis of
$\ell_1$-based methods.  One set of conditions, known as the
restricted isometry property~\cite{CandesTao06} or RIP for short, is
based on very strong constraints on the condition numbers of all
submatrices of $\Xmat$ up to size $2 \s$, requiring that they be
near-isometries (i.e., with condition numbers extremely close to $1$).
Such conditions are satisfied by matrices with columns that are all
very close to orthogonal (e.g., when $\Xmat$ has i.i.d. $N(0,1)$
entries and $\numobs = \Omega(\log {\pdim \choose 2 \spindex})$), but
are violated for many reasonable matrix classes (e.g., Toeplitz
matrices) that arise in statistical practice.  Zhang and
Huang~\cite{Huang06} imposed a weaker sparse Riesz condition, based on
imposing constraints (different from those of RIP) on the condition
numbers of all submatrices of $\Xmat$ up to a size that grows as a
function of $\s$ and $\numobs$.  Meinshausen and Yu~\cite{MeinYu06}
impose a bound in terms of the condition numbers or minimum and
maximum restricted eigenvalues for submatrices of $\Xmat$ up to size
$\s \log n$. It is unclear whether the conditions in Meinshausen and
Yu~\cite{MeinYu06} are weaker or stronger than the conditions in Zhang
and Huang~\cite{Huang06}.

The weakest known sufficient conditions to date are due to Bickel et
al.~\cite{BicRitTsy08}, who show that in addition to the column
normalization condition (Assumption~\ref{AssColNorm} in this paper),
it suffices to impose a milder condition, namely a lower bound on a
certain type of restricted eigenvalue (RE).  They show that this RE
condition is less restrictive than both the RIP
condition~\cite{CandesTao06} and the eigenvalue conditions imposed in
Meinshausen and Yu~\cite{MeinYu06}.  For a given vector $\mytheta \in
\real^\pdim$, let $\mytheta_{(j)}$ refer to the $j^{th}$ largest
coefficient in absolute value, so that we have the ordering
\begin{equation*}
\mytheta_{(1)} \; \geq \; \mytheta_{(2)} \; \geq \; \ldots \geq \;
\mytheta_{(\pdim-1)} \; \geq \; \mytheta_{(\pdim)}.
\end{equation*}
For a given scalar $c_0$ and integer $s = 1, 2, \ldots, \pdim$, let
define the set
\begin{eqnarray*}
\label{SubsetBickel}	
\Gamma(\s, c_0) & \defn & \biggr \{\mytheta \in \real^{\pdim} \, \mid
\; \sum_{j=\s + 1}^\pdim |\mytheta_{(j)}| \leq c_0 \sum_{j=1}^\s
|\mytheta_{(j)}| \biggr \}.
\end{eqnarray*}
In words, the set $\Gamma(\s, c_0)$ contains all vectors in
$\real^{\pdim}$ where the $\ell_1$-norm of the largest $\s$
co-ordinates provides an upper bound (up to constant $c_0$) to the
$\ell_1$ norm over the smallest $\pdim-\s$ co-ordinates.  For example
if $\pdim=3$, then the vector $(1,1/2, 1/4) \; \in \; \Gamma(1,1)$
whereas the vector $(1,3/4, 3/4)\notin \Gamma(1, 1)$.

With this notation, the restricted eigenvalue (RE) assumption can be
stated as follows:
\bass[Restricted lower eigenvalues~\cite{BicRitTsy08}]
\label{AssRE}
There exists a function $\kappa(\Xmat, c_0) > 0$ such that
\begin{eqnarray*}
\frac{1}{\sqrt{\numobs}} \| \Xmat \mytheta\|_2 & \geq & \kappa(\Xmat,
c_0)\| \mytheta \|_2 \quad \mbox{for all $\mytheta \in \Gamma(\s,
c_0)$.}
\end{eqnarray*}
\eass 
Bickel et. al~\cite{BicRitTsy08} require a slightly stronger condition
for bounding the $\ell_2$-loss in if $\s$ depends on $n$. However
the conditions are equivalent for fixed $\s$ and Assumption
\ref{AssRE} is much simpler to analyze and compare to Assumption
\ref{AssLower}. At this point, we have not seen conditions weaker than
Assumption \ref{AssRE}.

The following corollary of Proposition~\ref{PropRandDesign} shows that
Assumption~\ref{AssRE} is satisfied with high probability for broad
classes of Gaussian random designs:
\bcors
Suppose that $\rhomax(\CovMat)$ remains bounded, $\min_{v \in \BallZ(2
\s)} \frac{\|\CovMatSqrt v \|_2}{\|v\|_2} > 0$ and that 
\mbox{$\numobs
> \plaincon_3 \s \log \pdim$} for a sufficiently large constant. Then
a randomly drawn design matrix $\Xmat \in \real^{\numobs \times
\pdim}$ with i.i.d. $N(0, \CovMat)$ rows satisfies
Assumption~\ref{AssRE} with probability greater than $1 - \plaincon_1
\exp(-\plaincon_2 \numobs)$.
\ecors
\spro
Note that for any vector $\mytheta \in\Gamma(\s, c_0)$, we have
\begin{align*}
\|\theta\|_1 & \leq (1 + c_0) \sum_{j=1}^\s |\mytheta_{(j)}| \; \leq
\; (1 + c_0) \sqrt{\s} \|\mytheta\|_2.
\end{align*}
Consequently, if the bound~\eqref{EqnLowerFinal} holds, we have
\begin{align*}
\frac{\|X v\|_2}{\sqrt{\numobs}} & \geq \Big \{ \frac{\|\CovMatSqrt v
\|_2}{ 2 \|v\|_2} - 6 (1+ c_0) \big(\frac{\rhomax(\CovMat) 
\s \log
\pdim}{\numobs} \big)^{1/2} \Big \} \|v\|_2.
\end{align*}
Since we have assumed that \mbox{$\numobs > \plaincon_3 \s \log
\pdim$} for a sufficiently large constant,
the claim follows.
\fpro
Combined with the discussion following
Proposition~\ref{PropRandDesign}, this result shows that both the
conditions required by Theorem~\ref{ThmUpper} of this paper and the
analysis of Bickel et al.~\cite{BicRitTsy08} (both in the case $\qpar
= 0$) hold with high probability for Gaussian random designs.

\subsubsection{Comparison of RE assumption with Assumption~\ref{AssLower}}

In the case $\qpar = 0$, the condition required by the estimator that
performs least-squares over the $\ell_0$-ball---namely, the form of
Assumption~\ref{AssLower} used in Theorem~\ref{ThmUpper}(b)---is not
stronger than Assumption~\ref{AssRE}.  This fact was previously
established by Bickel et al. (see p.7,~\cite{BicRitTsy08}).  We now
provide a simple pedagogical example to show that the $\ell_1$-based
relaxation can fail to recover the true parameter while the optimal
$\ell_0$-based algorithm succeeds. In particular, let us assume that
the noise vector $w = 0$, and consider the design matrix
\begin{equation*}
\Xmat = \begin{bmatrix} 1 & -2 & -1 \\ 2 & -3 & -3 \end{bmatrix},
\end{equation*}
corresponding to a regression problem with $\numobs = 2$ and $\pdim =
3$.  Say that the regression vector $\betastar \in \real^3$ is hard
sparse with one non-zero entry (i.e., $\s = 1$).  Observe that the
vector $\Delta \defn \begin{bmatrix} 1 & 1/3 & 1/3
\end{bmatrix}$ belongs to the null-space of $\Xmat$, and moreover
$\Delta \in \Gamma(1,1)$ but $\Delta \notin \mathbb{B}_0(2)$. All the
$2 \times 2$ sub-matrices of $X$ have rank two, we have
$\mathbb{B}_0(2) \cap \Ker(X) = \{ 0 \}$, so that by known results
from Cohen et. al.~\cite{Cohen} (see, in particular, their Lemma 3.1),
the condition $\mathbb{B}_0(2) \cap \Ker(X) = \{ 0 \}$ implies that
the $\ell_0$-based algorithm can exactly recover any $1$-sparse
vector. On the other hand, suppose that, for instance, the true
regression vector is given by $\betaTr = \begin{bmatrix} 1 & 0 & 0
\end{bmatrix}$, If applied to this problem with no noise, the Lasso
would incorrectly recover the solution $\betahat \defn \begin{bmatrix}
0 & -1/3 & -1/3 \end{bmatrix}$ since $\|\betahat \|_1 = 2/3 \leq 1 =
\|\betaTr\|_1$.  Although this example is low-dimensional ($(\s,
\pdim) = (1,3)$), we suspect that higher dimensional examples of
design matrices that satisfy the conditions required for the minimax
rate but not satisfied for $\ell_1$-based methods may be constructed
using similar arguments. This construction highlights that there are
instances of design matrices $\Xmat$ for which $\ell_1$-based methods
fail to recover the true parameter $\betaTr$ for $\qpar = 0$ while the
optimal $\ell_0$-based algorithm succeeds.

In summary, for the hard sparsity case $\qpar = 0$, methods based on
$\ell_1$-relaxation can achieve the minimax rate $\order \big(\frac{\s
\log \pdim}{\numobs} \big)$ for $\ell_2$-error, but the current
analyses of these
$\ell_1$-methods~\cite{CandesTao06,MeinYu06,BicRitTsy08} are based on
imposing stronger conditions on the design matrix $\Xmat$ than those
required by the ``optimal'' estimator that performs least-squares over
the $\ell_0$-ball.

%%%%%%%%%%%%%%%%%%%%%%%%%%%%%%%%%%%%%%%%%%%%%%%%%%%%%%%%%%%%%%%%%%%%%%%%

\section{Proofs of main results}
\label{SecProof}

In this section, we provide the proofs of our main theorems, with more
technical lemmas and their proofs deferred to the appendices. To
begin, we provide a high-level overview that outlines the main steps
of the proofs.

\paragraph{Basic steps for lower bounds}
The proofs for the lower bounds follow an information-theoretic method
based on Fano's inequality~\cite{Cover}, as used in classical work on
nonparametric estimation~\cite{IbrHas81,YanBar99,Yu}.  A key
ingredient is a fine characterization of the metric entropy structure
of $\ell_\qpar$ balls~\cite{Kuhn,Carl90}.  At a high-level, the proof
of each lower bound follows the following three basic steps:

\begin{enumerate}
\item[(1)] Let $\|\cdot\|_\ast$ be the norm for which we wish to lower
bound the minimax risk; for Theorem~\ref{ThmLower}, the norm
$\|\cdot\|_\ast$ corresponds to the $\ell_{\rpar}$ norm, whereas for
Theorem~\ref{ThmLowerPred}, it is the $\ell_2$-prediction norm (the
square root of the prediction loss).  We first construct an
$\epspack$-packing set for $\Ballq(\myrad)$ in the norm
$\|\cdot\|_\ast$, where $\epspack > 0$ is a free parameter to be
determined in a later step.  The packing set is constructed by
deriving lower bounds on the packing numbers for $\Ballq(\myrad)$; we
discuss the concepts of packing sets and packing numbers at more
length in Section~\ref{MetricEntropy}. For the case of
$\ell_\qpar$-balls for $\qpar>0$, tight bounds on the packing numbers
in $\ell_{\rpar}$ norm have been developed in the approximation theory
literature~\cite{Kuhn}. For $\qpar = 0$, we use combinatorial to bound
the packing numbers. We use Assumption~\ref{AssLower} in order to
relate the packing number in the $\ell_2$-prediction norm to the
packing number in $\ell_2$-norm.   

\item [(2)] The next step is to use a standard reduction to show that
any estimator with minimax risk $\order(\epspack^2)$ must be able to
solve a hypothesis-testing problem over the packing set with vanishing
error probability.  More concretely, suppose that an adversary places
a uniform distribution over the $\epspack$-packing set in
$\Ballq(\myrad)$, and let this random variable be $\Theta$.  The
problem of recovering $\Theta$ is a multi-way hypothesis testing
problem, so that we may apply Fano's inequality to lower bound
the probability of error.  The Fano bound involves the log packing
number and the mutual information $I(\y; \Theta)$ between the observation
vector $y \in \real^\numobs$ and the random parameter $\Theta$ chosen
uniformly from the packing set.
\item[(3)] Finally, following a technique introduced by Yang and
Barron~\cite{YanBar99}, we derive an upper bound on the mutual
information between $\y$ and $\Theta$ by constructing an
$\epskey$-covering set for $\Ballq(\myrad)$ with respect to the
$\ell_2$-prediction semi-norm.  Using Lemma~\ref{LemKL} in
Section~\ref{MetEntrConvHull}, we establish a link between covering
numbers in $\ell_2$-prediction semi-norm to covering numbers in
$\ell_2$-norm.  Finally, we choose the free parameters $\epspack > 0$
and and $\epskey > 0 $ so as to optimize the lower bound.
\end{enumerate}

\paragraph{Basic steps for upper bounds}

The proofs for the upper bounds involve direct analysis of the natural
estimator that performs least-squares over the
$\ell_{\qpar}$-ball. The proof is constructive and involves two steps,
the first of which is standard while the second step is more specific
to the problem at hand:

\begin{enumerate}
\item[(1)] Since the estimator is based on minimizing the
least-squares loss over the ball $\Ballq(\Rq)$, some straightforward
algebra allows us to upper bound the $\ell_2$-prediction error by a
term that measures the supremum of a Gaussian empirical process over
the ball $\Ballq(2 \Rq)$. This step is completely generic and applies
to any least-squares estimator involving a linear model.
\item [(2)] The second and more challenging step involves computing
upper bounds on the supremum of the Gaussian process over $\Ballq(2
\Rq)$.  For each of the upper bounds, our approach is slightly
different in the details.  Common steps include upper bounds on the
covering numbers of the ball $\Ballq(2 \Rq)$, as well as on the image
of these balls under the mapping $\Xmat: \real^\pdim \rightarrow
\real^\numobs$.  For the case $\qpar = 1$, we make use of
Lemma~\ref{LemElloneElltwo} in order to relate the $\ell_1$-norm to
the $\ell_2$-norm for vectors that lie in an $\ell_{\qpar}$-ball.  For
$\qpar \in (0,1)$, we make use of some chaining and peeling results
from empirical process theory (e.g., Van de Geer ~\cite{vandeGeer}).
\end{enumerate}

\subsection{Packing, covering, and metric entropy}

\label{MetricEntropy}

The notion of packing and covering numbers play a crucial role in our
analysis, so we begin with some background, with emphasis on the
case of covering/packing for $\ell_\qpar$-balls.
\bdes[Covering and packing numbers]
Consider a metric space consisting of a set $\GenSet$ and a metric
$\rho: \GenSet \times \GenSet \rightarrow \real_+$.
\begin{enumerate}
\item[(a)] An $\epsilon$-covering of $\GenSet$ in the metric $\rho$ is
a collection $\{\beta^1, \ldots, \beta^\CovNum\} \subset \GenSet$ such
that for all $\beta \in S$, there exists some $i \in \{1, \ldots,
\CovNum\}$ with $\rho(\beta, \beta^i) \leq \epsilon$.  The
$\epsilon$-covering number $\CovNum(\epsilon; \GenSet, \rho)$ is the
cardinality of the smallest $\epsilon$-covering.
\item[(b)] A $\delta$-packing of $\GenSet$ in the metric $\rho$ is a
collection $\{\beta^1, \ldots, \beta^\PackNum \} \subset S$ such that
$\rho(\beta^i, \beta^j) \geq \delta$ for all $i \neq j$.  The
$\delta$-packing number $\PackNum(\delta; \GenSet, \rho)$ is the
cardinality of the largest $\delta$-packing.
\end{enumerate}
\edes
In simple terms, the covering number $\CovNum(\epsilon; \GenSet,
\rho)$ is the minimum number of balls with radius $\epsilon$ under the
metric $\rho$ required to completely cover the space, so that every
point in $\GenSet$ lies in some ball. The packing number
$\PackNum(\delta; \GenSet, \rho)$ is the maximum number of balls of
radius $\delta$ under metric $\rho$ that can be packed into the space
so that there is no overlap between any of the balls.  It is worth
noting that the covering and packing numbers are (up to constant
factors) essentially the same.  In particular, the inequalities
$\PackNum(\epsilon; \GenSet, \rho) \; \leq \; \CovNum(\epsilon;
\GenSet, \rho) \; \leq \; \PackNum(\epsilon/2; \GenSet, \rho)$ are
standard (e.g.,~\cite{Pollard84}).  Consequently, given upper and
lower bounds on the covering number, we can immediately infer similar
upper and lower bounds on the packing number.  Of interest in our
results is the logarithm of the covering number $\log
\CovNum(\epsilon; \GenSet, \rho)$, a quantity known as the
\emph{metric entropy}.

A related quantity, frequently used in the operator theory
literature~\cite{Kuhn,Sch84,Carl90}, are the (dyadic) entropy numbers
$\dyaeps_k(\GenSet; \rho)$, defined as follows for $k=1, 2, \ldots$
\begin{eqnarray}
\label{EqnDefnDyaEps}
\dyaeps_k(\GenSet; \rho) & = & \inf \big \{ \epsilon > 0 \, \mid \,
\CovNum(\epsilon; \GenSet, \rho) \leq 2^{k-1} \big \}.
\end{eqnarray}
By definition, note that we have $\dyaeps_k(\GenSet; \rho) \leq
\delta$ if and only if $\log N(\delta; \GenSet, \rho) \leq k$.

\subsubsection{Metric entropies of $\ell_\qpar$-balls}

Central to our proofs is the metric entropy of the ball
$\Ballq(\myrad)$ when the metric $\rho$ is the $\ell_\rpar$-norm, a
quantity which we denote by $\log \CovNum_{\rpar, \qpar}(\epsilon)$.
The following result, which provides upper and lower bounds on this
metric entropy that are tight up to constant factors, is an adaptation
of results from the operator theory literature~\cite{Kuhn,GueLit00};
see Appendix~\ref{AppLemEntropyBounds} for the details.  All bounds
stated here apply to a dimension $\pdim \geq 2$.
\blems
\label{LemEntropyBounds} 
Assume that $\qpar \in (0, 1]$ and $\rpar \in [1, \infty]$ with $\rpar
> \qpar$.  Then there is a constant $\UqrOne$, depending only on
$\qpar$ and $\rpar$, such that
\begin{equation}
\label{EqnEntropyBoundsUpper}
\log \CovNum_{\rpar, \qpar} (\epsilon) \; \leq \; \UqrOne \, \Big[
\Rq^{\frac{\rpar}{\rpar-\qpar}} \big(\frac{1}{\epsilon}
\big)^{\frac{\rpar \qpar}{\rpar-\qpar}} \log \pdim \Big] \quad
\mbox{for all $\epsilon \in (0, \Rq^{1/\qpar})$.}
\end{equation}
Conversely, suppose in addition that $\epsilon < 1$ and
$\epsilon^\rpar = \Omega\big(\frac{\log \pdim}{\pdim^\kOne}
\big)^{\frac{\rpar- \qpar}{\qpar}}$ for some fixed $\kOne \in (0,1)$,
depending only on $\qpar$ and $\rpar$.  Then there is a constant
$\LqrOne \leq \UqrOne$, depending only on $\qpar$ and $\rpar$, such
that
\begin{equation}
\label{EqnEntropyBoundsLower}
\log \CovNum_{\rpar, \qpar} (\epsilon) \; \geq \; \LqrOne \;
\Big[\Rq^{\frac{\rpar}{\rpar-\qpar}} \big(\frac{1}{\epsilon}
\big)^{\frac{\rpar \qpar}{\rpar-\qpar}} \log \pdim \Big].
\end{equation}
\elems
\newcommand{\REFLEMENTROPY}{3}

\noindent {\bf{Remark:}} In our application of the lower
bound~\eqref{EqnEntropyBoundsLower}, our typical choice of
$\epsilon^\rpar$ will be of the order $\order\big(\frac{\log
\pdim}{\numobs} \big)^{\frac{\rpar- \qpar}{2}}$. It can be verified
that as long as there exists a $\kappa \in (0,1)$ such that $\frac{\pdim}{\Rq \numobs^{\qpar/2}} = \Omega (\pdim^\kappa)$ (which is stated at the beginning of Section~\ref{SecMain}) and $\rpar > \qpar$, then there exists some fixed $\kOne \in (0,1)$, depending only on $\rpar$ and $\qpar$,
such that $\epsilon$ lies in the range required for the lower
bound~\eqref{EqnEntropyBoundsLower} to be valid.

\subsubsection{Metric entropy of $\qpar$-convex hulls}
\label{MetEntrConvHull}

The proofs of the lower bounds all involve the Kullback-Leibler (KL)
divergence between the distributions induced by different parameters
$\beta$ and $\beta'$ in $\Ballq(\myrad)$.  Here we show that for the
linear observation model~\eqref{EqnLinearObs}, these KL divergences
can be represented as $\qpar$-convex hulls of the columns of the
design matrix, and provide some bounds on the associated metric
entropy.

For two distributions $\mprob$ and $\qprob$ that have densities $d
\mprob$ and $d \qprob$ with respect to some base measure $\mu$, the
Kullback-Leibler (KL) divergence is given by $\kull{\mprob}{\qprob} =
\int \log \frac{d \mprob}{d \qprob} \; \mprob(d \mu)$.  We use
$\mprob_{\beta}$ to denote the distribution of $y \in \real$ under the
linear regression model---in particular, it corresponds to the
distribution of a $N(\Xmat \beta, \sigma^2 I_{\numobs \times
\numobs})$ random vector. A straightforward computation then leads to
\begin{eqnarray}
\kull{\mprob_{\beta}}{\mprob_{\beta'}} & = & \frac{1}{2 \sigma^2} \,
\|\Xmat \beta - \Xmat \beta'\|_2^2.
\end{eqnarray}

Therefore, control of KL-divergences requires understanding of the
metric entropy of the $\qpar$-convex hull of the rescaled columns of
the design matrix $\Xmat$---in particular, the set
\begin{eqnarray}
\label{EqnQparConv}
\absconv_\qpar(\Xmat/\sqrt{\numobs}) & \defn & \big \{
\frac{1}{\sqrt{\numobs}} \sum_{j=1}^\pdim \mytheta_j \Xmat_j\, \mid
\, \mytheta \in \Ballq(1) \big \}.
\end{eqnarray}
We have introduced the normalization by $1/\sqrt{\numobs}$ for later
technical convenience.

Under the column normalization condition, it turns out that the metric
entropy of this set with respect to the $\ell_2$-norm is essentially
no larger than the metric entropy of $\Ballq(\myrad)$, as summarized
in the following
\blems
\label{LemKL}
Suppose that $\Xmat$ satisfies the column normalization condition
(Assumption~\ref{AssColNorm} with constant $\colnorm$).  Then there is
a constant $U'_{\qpar,2}$ depending only on $\qpar \in (0,1]$ such
that
\begin{eqnarray*}
\log \CovNum(\epsilon, \absconv_\qpar(\Xmat/\sqrt{\numobs}),
\|\cdot\|_2) & \leq & U'_{\qpar, 2} \; \Big[ \Rq^{\frac{2}{2-\qpar}}
\big(\frac{\colnorm}{\epsilon} \big)^{\frac{2 \qpar}{2-\qpar}} \log
\pdim \Big].
\end{eqnarray*}
\elems
\newcommand{\REFLEMKL}{4}
\noindent The proof of this claim is provided in
Appendix~\ref{AppLemKL}.  Note that apart from a different constant,
this upper bound on the metric entropy is identical to that for $\log
N_{2, \qpar}(\epsilon/\colnorm)$ from Lemma~\ref{LemEntropyBounds}.
Up to constant factors, this upper bound cannot be tightened in
general (e.g., consider $\numobs = \pdim$ and $\Xmat = I$).

\subsection{Proof of lower bounds}

We begin by proving our main results that provide lower bounds on
minimax risks, namely Theorems~\ref{ThmLower} and~\ref{ThmLowerPred}.

\subsubsection{Proof of Theorem~\ref{ThmLower}}
\label{SecProofThmLower}

Recall that the lower bounds in Theorem~\ref{ThmLower} are the maximum
of two expressions, one corresponding to the diameter of the set $\Kx$
intersected with the $\ell_\qpar$-ball, and the other correspond to
the metric entropy of the $\ell_\qpar$-ball.

We begin by deriving the lower bound based on the diameter of $\Kx =
\Ballq(\myrad) \cap \Ker(\Xmat)$.  The minimax risk is lower bounded
as
\begin{eqnarray*}
\min_{\betahat} \max_{\beta \in \Bq} \Exs \|\betahat -
\beta\|_{\rpar}^{\rpar} & \geq & \min_{\betahat} \max_{\beta \in \Kx}
\Exs \|\betahat - \beta\|_{\rpar}^{\rpar},
\end{eqnarray*}
where the inequality follows from the inclusion $\Kx \subseteq
\Ballq(\myrad)$.  For any $\beta \in \Kx$, we have $\y = \Xmat \beta +
\w = \w$, so that $\y$ contains no information about $\beta \in \Kx$.
Consequently, once $\betahat$ is chosen, the adversary can always
choose an element $\beta \in \Kx$ such that \mbox{$\|\betahat -
\beta\|_{\rpar} \geq \frac{1}{2} \diam_{\rpar}(\Kx)$.}  Indeed, if
$\|\betahat\|_\rpar \geq \frac{1}{2} \diam_{\rpar}(\Kx)$, then the
adversary chooses $\beta = 0 \in \Kx$.  On the other hand, if
$\|\betahat \|_\rpar \leq \frac{1}{2} \diam_\rpar(\Kx)$, then the
adversary can choose some $\beta \in \Kx$ such that $\|\beta\|_\rpar =
\diam_\rpar(\Kx)$.  By triangle inequality, we then have $\|\beta -
\betahat \|_\rpar \geq \|\beta\|_\rpar - \|\betahat\|_\rpar \geq
\frac{1}{2} \diam_\rpar(\Kx)$.  Overall, we conclude that
\begin{eqnarray*}
\min_{\betahat} \max_{\beta \in \Bq} \Exs \|\betahat -
\beta\|_{\rpar}^{\rpar} & \geq & \big(\frac{1}{2}
\diam_{\rpar}(\Kx)\big)^{\rpar}.
\end{eqnarray*}
In the following subsections, we establish the second terms in the
lower bounds via the Fano method, a standard approach for minimax lower
bounds. Our proofs of part (a) and (b) are based on slightly
different arguments.

\myparagraph{Proof of Theorem~\ref{ThmLower}(a)}

Let $\PackNum = \PackNum_\rpar(\epspack)$ be the cardinality of a
maximal packing of the ball $\Ballq(\myrad)$ in the $\ell_\rpar$
metric, say with elements $\{ \beta^1, \ldots, \beta^\PackNum \}$.  A
standard argument (e.g.,~\cite{Hasminskii,YanBar99,Yu}) yields a lower
bound on the minimax $\ell_\rpar$-risk in terms of the error in a
multi-way hypothesis testing problem: in particular, we have
\begin{eqnarray*}
\min_{\betahat} \max_{\beta \in \Bq} \Exs \|\betahat -
\beta\|_\rpar^\rpar \geq \frac{1}{2^\rpar} \; \epspack^{\rpar} \;
\min_{\betatil} \mprob[\betatil \neq \betaRand]
\end{eqnarray*}
where the random vector $\betaRand \in \real^\pdim$ is uniformly
distributed over the packing set $\{\beta^1, \ldots,
\beta^\PackNum\}$, and the estimator $\betatil$ takes values in the
packing set.  Applying Fano's inequality~\cite{Cover} yields the lower
bound
\begin{eqnarray}
\label{EqnFano} 
\mprob[\betaRand \neq \betatil] \geq 1 - \frac{I(\betaRand; \y)+\log
2}{\log \PackNum_\rpar(\epspack)},
\end{eqnarray}
where $I(\betaRand; \y)$ is the mutual information between random
parameter $\betaRand$ in the packing set and the observation vector
$\y \in \real^\numobs$.

It remains to upper bound the mutual information; we do so by
following the procedure of Yang and Barron~\cite{YanBar99}, which is
based on covering the model space $\{ \mprob_{\beta}, \; \beta \in
\Ballq(\myrad) \}$ under the square-root Kullback-Leibler divergence.
As noted prior to Lemma~\ref{LemKL}, for the Gaussian models given
here, this square-root KL divergence takes the form
$\sqrt{\kull{\mprob_\beta}{\mprob_{\beta'}}} = \frac{1}{\sqrt{2
\sigma^2}} \| \Xmat(\beta - \beta')\|_2$.  Let $\CovNum =
\CovNum_2(\epskey)$ be the minimal cardinality of an $\epskey$-covering
of $\Ballq(\myrad)$ in $\ell_2$-norm.  Using the upper bound on the
dyadic entropy of $\absconv_\qpar(\Xmat)$ provided by
Lemma~\ref{LemKL}, we conclude that there exists a set $\{ \Xmat
\beta^1, \ldots, \Xmat \beta^\CovNum \}$ such that for all $\Xmat
\beta \in \absconv_\qpar(\Xmat)$, there exists some index $i$ such
that $\|\Xmat(\beta - \beta^i)\|_2/\sqrt{\numobs} \leq \gencon \,
\colnorm \, \epskey$.  Following the argument of Yang and
Barron~\cite{YanBar99}, we obtain that the mutual information is upper
bounded as
\begin{eqnarray*}
\label{EqnInfoUpper}
I(\betaRand; \y) & \leq & \log \CovNum(\epskey) + \frac{\gencon^2 \,
\numobs}{\sigma^2} \colnormsq \epskey^2.
\end{eqnarray*}
Combining this upper bound with the Fano lower bound~\eqref{EqnFano}
yields
\begin{eqnarray}
\label{EqnMatchProb}
\Prob[\betaRand \neq \betatil] & \geq & 1 - \frac{\log
\CovNum_2(\epskey) + \frac{\gencon^2 \, \numobs}{\sigma^2} \colnormsq
\, \epskey^2 + \log 2}{\log \PackNum_\rpar(\epspack)}.
\end{eqnarray}
The final step is to choose the packing and covering radii ($\epspack$
and $\epskey$ respectively) such that the lower
bound~\eqref{EqnMatchProb} remains strictly above zero, say bounded
below by $1/4$.  In order to do so, suppose that we choose the pair
$(\epskey, \epspack)$ such that
\begin{subequations}
\label{EqnKeyRelation}
\begin{eqnarray}
\label{EqnKeyRelationA}
\frac{\gencon^2 \, \numobs}{\sigma^2} \colnormsq \, \epskey^2
& \leq & \log \CovNum_2(\epskey), \quad \mbox{and} \\
\label{EqnKeyRelationB}
\log \PackNum_\rpar(\epspack) & \leq  &  4 \log
\CovNum_2(\epskey).
\end{eqnarray}
\end{subequations}
As long as $\CovNum_2(\epskey) \geq 2$, we are then guaranteed that
\begin{eqnarray}
\label{EqnOneFourth}
\Prob[\betaRand \neq \betatil] & \geq & 1 - \frac{2 \log
\CovNum_2(\epskey) + \log 2}{4 \log \CovNum_2(\epskey)} \; \geq \;
1/4,
\end{eqnarray}
as desired.

It remains to determine choices of $\epskey$ and $\epspack$ that
satisfy the relations~\eqref{EqnKeyRelation}.  From
Lemma~\ref{LemEntropyBounds}, relation~\eqref{EqnKeyRelationA} is
satisfied by choosing $\epskey$ such that $\frac{ \gencon^2 \,
\numobs}{\sigma^2} \colnormsq \, \epskey^2 \, = \, \MyLqTwo \;
\Big[\Rq^{\frac{2}{2-\qpar}} \big(\frac{1}{\epskey} \big)^{\frac{2
\qpar}{2 - \qpar}} \log \pdim \Big]$, or equivalently such that
\begin{eqnarray*}
\big( \epskey\big)^{\frac{4}{2-\qpar}} & = & \Theta \big(
\Rq^{\frac{2}{2-\qpar}} \frac{\sigma^2}{\colnormsq} \; \frac{\log
\pdim}{\numobs} \big).
\end{eqnarray*}

In order to satisfy the bound~\eqref{EqnKeyRelationB}, it suffices
to choose $\epspack$ such that
\begin{eqnarray*}
\MyUq \; \Big[\Rq^{\frac{\rpar}{\rpar-\qpar}} \big(\frac{1}{\epspack}
\big)^{\frac{\rpar \qpar}{\rpar - \qpar}} \log \pdim \Big] & \leq & 4
\MyLqTwo \; \Big[\Rq^{\frac{2}{2-\qpar}} \big(\frac{1}{\epskey}
\big)^{\frac{2 \qpar}{2 - \qpar}} \log \pdim \Big],
\end{eqnarray*}
or equivalently such that
\begin{eqnarray*}
\epspack^\rpar & \geq & \big[ \frac{\MyUq}{4 \MyLqTwo}
\big]^{\frac{\rpar-\qpar}{\qpar}} \; \biggr \{ \big(\epskey
\big)^{\frac{4}{2-\qpar}} \biggr \}^{\frac{\rpar - \qpar}{2}} \Rq^{\frac{2-\rpar}{2-\qpar}}\\
%
%& = &  \big[ \frac{\MyUq}{4 \MyLqTwo}
%\big]^{\frac{\rpar-\qpar}{\qpar}} \;  \Big[\MyLqTwo \;
%\Rq^{\frac{2}{2-\qpar}} \frac{\sigma^2}{\colnormsq} \;
%\frac{\log \pdim}{\numobs} \Big]^{\frac{\rpar-\qpar}{2}} \\
%
& = & \big[ \frac{\MyUq}{4 \MyLqTwo} \big]^{\frac{\rpar-\qpar}{\qpar}}
\; \MyLqTwo^{\frac{\rpar-\qpar}{2}} \; \Rq \; \Big[\frac{\sigma^2}{\colnormsq} \;
\frac{\log \pdim}{\numobs} \Big]^{\frac{\rpar-\qpar}{2}}
\end{eqnarray*}
Combining this bound with the lower bound~\eqref{EqnOneFourth} on the
hypothesis testing error probability and substituting into
equation~\eqref{EqnLowerBound}, we obtain
\begin{eqnarray*}
\min_{\betahat} \max_{\beta \in \Bq} \Exs \|\betahat -
\beta\|_\rpar^\rpar & \geq & \AnnoyLower \;
\Rq \;
\Big[\frac{\sigma^2}{\colnormsq} \; \frac{\log \pdim}{\numobs}
\Big]^{\frac{\rpar-\qpar}{2}},
\end{eqnarray*}
which completes the proof of Theorem~\ref{ThmLower}(a).

\myparagraph{Proof of Theorem~\ref{ThmLower}(b)}

In order to prove Theorem~\ref{ThmLower}(b), we require some definitions and an auxiliary lemma. For any integer $\spindex \in \{1, \ldots, \pdim \}$,
we define the set
\begin{eqnarray*}
\Hyper(\spindex) & \defn & \big \{ z \in \{-1, 0, +1 \}^\pdim \, \mid
\, \|z\|_0 = \spindex \big \}.
\end{eqnarray*}
Although the set $\Hyper$ depends on $\spindex$, we frequently
drop this dependence so as to simplify notation.  We define the
Hamming distance \mbox{$\rhohamm(z, z') = \sum_{j=1}^\pdim \Ind[z_j
\neq z'_j]$} between the vectors $z$ and $z'$.  We prove the following
result in Appendix~\ref{AppLemHyperCube}:
\blems
\label{LemHyperCube}
There exists a subset $\HyperSub \subset \Hyper$ with cardinality
$|\HyperSub| \geq \exp(\frac{\spindex}{2} \log \frac{\pdim -
\spindex}{\spindex/2})$ such that $\rhohamm(z, z') \geq
\frac{\spindex}{2}$ for all $z, z' \in \HyperSub$.
\elems
\newcommand{\REFLEMHYPER}{5}
Now consider a rescaled version of the set $\HyperSub$, say
$\sqrt{\frac{2}{\spindex}} \epspack \HyperSub$ for some $\epspack > 0$
to be chosen.  For any elements $\beta, \beta' \in
\frac{\epspack}{\sqrt{\spindex}} \HyperSub$, we have the following
bounds on the $\ell_2$-norm of their difference:
\begin{subequations}
\begin{eqnarray}
\label{EqnLowerBeta}
\| \beta - \beta'\|^2_2 & \geq & \epspack^2, \quad \mbox{and} \\
\label{EqnUpperBeta}
\| \beta - \beta'\|^2_2 & \leq & 8 \epspack^2.
\end{eqnarray}
\end{subequations}
Consequently, the rescaled set $\sqrt{\frac{2}{\spindex}} \epspack
\HyperSub$ is an $\epspack$-packing set in $\ell_2$ norm with
$\PackNum_2(\epspack) = |\HyperSub|$ elements, say $\{\beta^1, \ldots,
\beta^\PackNum\}$. Using this packing set, we now follow the same
classical steps as in the proof of Theorem~\ref{ThmLower}(a), up until
the Fano lower bound~\eqref{EqnFano}.

At this point, we use an alternative upper bound on the mutual
information, namely the bound $I(\y; \betaRand) \leq
\frac{1}{{\PackNum \choose 2}} \sum_{i \neq j}
\kull{\beta^i}{\beta^j}$, which follows from the convexity of mutual
information~\cite{Cover}.  For the linear observation
model~\eqref{EqnLinearObs}, we have $\kull{\beta^i}{\beta^j} =
\frac{1}{2 \sigma^2} \|\Xmat(\beta^i - \beta^j)\|_2^2$.  Since $(\beta
- \beta') \in \Ball_0(2 \s)$ by construction, from the assumptions on
$\Xmat$ and the upper bound bound~\eqref{EqnUpperBeta}, we conclude
that
\begin{eqnarray*}
I(\y; \betaRand) & \leq & \frac{8 \numobs \conuppersq \, \epspack^2}{2
\sigma^2}.
\end{eqnarray*}
Substituting this upper bound into the Fano lower
bound~\eqref{EqnFano}, we obtain
\begin{eqnarray*}
\mprob[\betaRand \neq \betatil] & \geq & 1 - \frac{\frac{ 8 \, \numobs
\conuppersq}{2 \sigma^2} \epspack^2 + \log(2)}{ \frac{\spindex}{2}
\log \frac{\pdim - \spindex}{\spindex/2}}.
\end{eqnarray*}
Setting $\epspack^2 = \frac{1}{32} \frac{\sigma^2}{\conuppersq}
\frac{\spindex}{2 n} \log \frac{ \pdim - \spindex}{\spindex/2}$ ensures
that this probability is at least $1/4$.  Consequently, combined with
the lower bound~\eqref{EqnLowerBound}, we conclude that
\begin{eqnarray*}
\min_{\betahat} \max_{\beta \in \Bq} \Exs \|\betahat -
\beta\|_\rpar^\rpar & \geq & \frac{1}{2^\rpar} \frac{1}{4}
(\frac{1}{32})^{\rpar/2} \; \Big[ \frac{\sigma^2}{\conuppersq}
\frac{\spindex}{2 n} \log \frac{ \pdim - \spindex}{\spindex/2}
\Big]^{\frac{\rpar}{2}}.
\end{eqnarray*}
As long as the ratio $\pdim/\spindex \geq 1 + \delta$ for some $\delta
> 0$ we have $\log(\pdim/\spindex - 1) \geq c \log(\pdim/\spindex)$
for some constant $c > 0$, from which the result follows.

%%%%%%%%%%%%%%%%%%%%%%%%%%%%%%%%%%%%%%%%%%%%%%%%%%%%%%%%%%%%%%%%%%%%%%%%%%

\subsubsection{Proof of Theorem~\ref{ThmLowerPred}}

We use arguments similar to the proof of Theorem~\ref{ThmLower}
in order to establish lower bounds on prediction error $\|\Xmat
(\betahat - \betastar)\|_2/\sqrt{\numobs}$.  

\myparagraph{Proof of Theorem~\ref{ThmLowerPred}(a)}
For some $\epspack^2 = \Omega(\myrad \, (\frac{\log
\pdim}{\numobs})^{1 - \qpar/2})$, let $\{ \beta^1, \ldots,
\beta^\PackNum \}$ be an $\epspack$ packing of the ball
$\Ballq(\myrad)$ in the $\ell_2$ metric, say with a total of $\PackNum
= \PackNum(\epspack/\colnorm)$ elements.  We first show that if
$\numobs$ is sufficiently large, then this set is also a $\conlower
\epspack/2$-packing set in the prediction (semi)-norm. From
Assumption~\ref{AssLower}, for each $i \neq j$,
\begin{eqnarray}
\label{EqnInitialLower}
\frac{\|\Xmat (\beta^i - \beta^j)\|_2}{\sqrt{\numobs}} & \geq &
\conlower \, \|\beta^i- \beta^j\|_2 - 
\SpecFunLower.
\end{eqnarray}

Using the assumed lower bound on $\epspack^2$---namely, $\epspack^2 =
\Omega\big(\Rq (\frac{\log \pdim}{n})^{1-\frac{\qpar}{2}})$---and the
initial lower bound~\eqref{EqnInitialLower}, we conclude that
\mbox{$\frac{\|\Xmat (\beta^i - \beta^j)\|_2}{\sqrt{\numobs}} \geq
\conlower \epspack/2$} once $\numobs$ is larger than some finite
number.

We have thus constructed a $\conlower \epspack/2$-packing set
in the (semi)-norm $\|\Xmat(\beta^i - \beta^j)\|_2$.  As in the proof
of Theorem~\ref{ThmUpper}(a), we follow a standard approach to reduce
the problem of lower bounding the minimax error to the error
probability of a multi-way hypothesis testing problem.  After this
step, we apply the Fano inequality to lower bound this error
probability via
\begin{eqnarray*}
\mprob[X \betaRand \neq X \betatil] \geq 1 - \frac{I(X \betaRand^i;
\y)+\log 2}{\log \PackNum_2(\epspack)},
\end{eqnarray*}
where $I(X \betaRand^i; \y)$ now represents the mutual
information\footnote{Despite the difference in notation, this mutual
information is the same as $I(\betaRand; \y)$, since it measures the
information between the observation vector $y$ and the discrete index
$i$.}  between random parameter $X \betaRand$ (uniformly distributed
over the packing set) and the observation vector $\y \in
\real^\numobs$.

From Lemma~\ref{LemKL}, the $\colnorm \, \epsilon$-covering number of
the set $\absconv_\qpar(\Xmat)$ is upper bounded (up to a constant
factor) by the $\epsilon$ covering number of $\Ballq(\myrad)$ in
$\ell_2$-norm, which we denote by $\CovNum_2(\epskey)$.  Following the
same reasoning as in Theorem~\ref{ThmUpper}(a), the mutual information
is upper bounded as
\begin{eqnarray*}
I(\Xmat \betaRand; \y) & \leq & \log \CovNum_2(\epskey) +
\frac{\numobs}{2 \sigma^2} \colnormsq \, \epskey^2.
\end{eqnarray*}
Combined with the Fano lower bound, we obtain
\begin{eqnarray}
\label{EqnMatchProbTwo}
\Prob[\Xmat \betaRand \neq \Xmat \betatil] & \geq & 1 - \frac{\log
\CovNum_2(\epskey) + \frac{ \numobs}{\sigma^2} \, \colnormsq \,
\epskey^2 + \log 2}{\log \PackNum_\rpar(\epspack)}.
\end{eqnarray}
Lastly, we choose the packing and covering radii ($\epspack$ and
$\epskey$ respectively) such that the lower
bound~\eqref{EqnMatchProbTwo} remains strictly above zero, say bounded
below by $1/4$.  It suffices to choose the pair $(\epskey, \epspack)$
to satisfy the relations~\eqref{EqnKeyRelationA}
and~\eqref{EqnKeyRelationB}.  As long as $\epskey^2 > \frac{\log
\pdim}{\numobs}$ and $\CovNum_2(\epskey) \geq 2$, we are then
guaranteed that
\begin{eqnarray*}
\Prob[\Xmat \betaRand \neq \Xmat \betatil] & \geq & 1 - \frac{2 \log
\CovNum_2(\epskey) + \log 2}{4 \log \CovNum_2(\epskey)} \; \geq \;
1/4,
\end{eqnarray*}
as desired.  Recalling that we have constructed a $\epspack
\conlower/2$ covering in the prediction (semi)-norm, we obtain
\begin{eqnarray*}
\min_{\betahat} \max_{\beta \in \Bq} \Exs \|\Xmat \, (\betahat -
\beta) \|_2^2/n & \geq & \Annoy'_{2, \qpar} \; \Rq \; \conlowersq \;
\Big[\frac{\sigma^2}{\colnormsq} \; \frac{\log \pdim}{\numobs}
\Big]^{1 - \qpar/2},
\end{eqnarray*}
for some constant $\Annoy'_{2, \qpar} > 0$. This completes the proof of
Theorem~\ref{ThmLowerPred}(a).

\myparagraph{Proof of Theorem~\ref{ThmLowerPred}(b)}

Recall the assertion of Lemma~\ref{LemHyperCube}, which guarantees the
existence of a set $\frac{\epspack^2}{2 \spindex} \HyperSub$ is an
$\epspack$-packing set in $\ell_2$-norm with $\PackNum_\rpar(\epspack)
= |\HyperSub|$ elements, say $\{\beta^1, \ldots, \beta^\PackNum\}$,
such that the bounds~\eqref{EqnLowerBeta} and~\eqref{EqnUpperBeta}
hold, and such that $\log |\HyperSub| \geq \frac{\spindex}{2} \log
\frac{\pdim - \spindex}{\spindex/2}$.  By construction, the difference
vectors $(\beta^i - \beta^j) \in \Ball_0(2 \spindex)$, so that by
assumption, we have
\begin{equation}
\label{EqnUpperXbeta}
\|\Xmat(\beta^i - \beta^j)\|/\sqrt{\numobs} \; \leq \; \conupper
\|\beta^i - \beta^j\|_2 \; \leq \; \conupper \sqrt{8} \; \epspack.
\end{equation}
In the reverse direction, since Assumption~\ref{AssLower} holds
with $\SpecFunLower = 0$, we have
\begin{eqnarray}
\label{EqnLowerXbeta}
\|\Xmat (\beta^i - \beta^j)\|_2/\sqrt{\numobs} & \geq & \conlower
\epspack.
\end{eqnarray}
We can follow the same steps as in the proof of
Theorem~\ref{ThmLower}(b), thereby obtaining an upper bound the mutual
information of the form $I( \Xmat \betaRand; y) \leq 8 \conuppersq
\numobs \epspack^2$.  Combined with the Fano lower bound, we have
\begin{eqnarray*}
\mprob[\Xmat \betaRand \neq \Xmat \betatil] & \geq & 1 - \frac{\frac{
8 \, \numobs \conuppersq}{2 \sigma^2} \epspack^2 + \log(2)}{
\frac{\spindex}{2 \numobs} \log \frac{\pdim - \spindex}{\spindex/2}}.
\end{eqnarray*}
Remembering the extra factor of $\conlower$ from the lower
bound~\eqref{EqnLowerXbeta}, we obtain the lower bound
\begin{eqnarray*}
\min_{\betahat} \max_{\beta \in \BZ} \Exs \frac{1}{\numobs}\|\Xmat(\betahat -
\beta)\|_2^2 & \geq & \Annoy'_{0, \qpar} \; \conlowersq \;
\frac{\sigma^2}{\conuppersq} \; \spindex \log \frac{ \pdim -
\spindex}{\spindex/2}.
\end{eqnarray*}
Repeating the argument from the proof of Theorem~\ref{ThmLower}(b)
allows us to further lower bound this quantity in terms of
$\log(\pdim/\spindex)$, leading to the claimed form of the bound.

%%%%%%%%%%%%%%%%%%%%%%%%%%%%%%%%%%%%%%%%%%%%%%%%%%%%%%%%%%%%%%%%%%%%%%%%%%%%%%
%%%%%%%%%%%%%%%%%%%%%%%%%%%%%%%%%%%%%%%%%%%%%%%%%%%%%%%%%%%%%%%%%%%%%%%%%%%%

\subsection{Proof of achievability results}

We now turn to the proofs of our main achievability results, namely
Theorems~\ref{ThmUpper} and~\ref{ThmUpperPred}, that provide upper
bounds on minimax risks. We prove all parts of these theorems by
analyzing the family of $M$-estimators
\begin{eqnarray*}
\betahat & \in & \arg \min_{\|\beta\|_\qpar^\qpar \leq \myrad} \|\y -
\Xmat \beta\|_2^2.
\end{eqnarray*}

We begin by deriving an elementary inequality that is useful
throughout the analysis.  Since the vector $\betastar$ satisfies the constraint $\|\betastar\|_\qpar^\qpar \leq \myrad$ meaning $\betastar$ is a feasible point, we have $\|\y - \Xmat \beta\|_2^2 \leq \|\y - \Xmat \betastar\|_2^2$.  Defining $\DeltaHat = \betahat - \betastar$ and
performing some algebra, we obtain the inequality
\begin{eqnarray}
\label{EqnBasic}
\frac{1}{\numobs} \|\Xmat \DeltaHat\|_2^2 & \leq & \frac{2 |\wnoise^T
  \Xmat \DeltaHat|}{\numobs}.
\end{eqnarray}

\subsubsection{Proof of Theorem~\ref{ThmUpper}}
\label{SecProofThmUpper}

We begin with the proof of Theorem~\ref{ThmUpper}, in which we upper
bound the minimax risk in squared $\ell_2$-norm.

\myparagraph{Proof of Theorem~\ref{ThmUpper}(a)}
To begin, we may apply Assumption~\ref{AssLower} to the inequality~\eqref{EqnBasic} to obtain
\begin{eqnarray*}
\big[\max ( 0, \; \conlower \|\DeltaHat\|_2 - \SpecFunLower) \big]^2 & \leq & 2 |\wnoise^T \Xmat
\DeltaHat|/\numobs \\
& \leq & \frac{2}{\numobs} \| \wnoise^T \Xmat\|_\infty
\|\DeltaHat\|_1.
\end{eqnarray*}
Since $w_i \sim N(0,\sigma^2)$ and the columns of $\Xmat$ are
normalized, each entry of $\frac{2}{\numobs} \wnoise^T \Xmat$ is
zero-mean Gaussian with variance at most $4 \sigma^2
\colnormsq/\numobs$.  Therefore, by union bound and standard Gaussian
tail bounds, we obtain that the inequality
\begin{eqnarray}
\label{EqnFirstQuad}
\big[\max ( 0, \; \conlower \|\DeltaHat\|_2 - \SpecFunLower) \big]^2 & \leq & 2 \sigma \colnorm \sqrt{\frac{3 \log
\pdim}{\numobs}} \, \|\DeltaHat\|_1
\end{eqnarray}
holds with probability greater than $1-\gencon_1 \exp(-\gencon_2
\numobs)$. Consequently, we may conclude that at least one of the two
following alternatives must hold
\begin{subequations}
\begin{eqnarray}
\label{EqnAltA}
\|\DeltaHat\|_2 & \leq & \frac{2\SpecFunLower}{\conlower}, \quad \mbox{or} \\
\label{EqnAltB}
\|\DeltaHat\|^2_2 & \leq & \frac{2 \sigma \colnorm}{\conlowersq}
\sqrt{\frac{3 \log \pdim}{\numobs}} \|\DeltaHat\|_1.
\end{eqnarray}
\end{subequations}

Suppose first that alternative~\eqref{EqnAltA} holds.  Consequently for we have
\begin{align*}
\|\DeltaHat\|_2^2 & \leq o\biggr( \myrad \big( \frac{\log
\pdim}{\numobs} \big)^{1 - \qpar/2} \biggr),
\end{align*}
which is the same up to constant rate than claimed in Theorem~\ref{ThmUpper}(a).

On the other hand, suppose that alternative~\eqref{EqnAltB} holds. Since both $\betahat$ and $\betastar$ belong to $\Ballq(\myrad)$, we
have $\|\DeltaHat\|_\qpar^\qpar  = \sum_{j=1}^{\pdim}{|\DeltaHat_j|^\qpar} \leq 2 \myrad$. Therefore we can exploit Lemma~\ref{LemElloneElltwo} by setting $\tau =
\frac{2 \sigma \colnorm}{\conlowersq} \sqrt{\frac{3 \log
\pdim}{\numobs}}$, thereby obtaining the bound $\|\DeltaHat\|_2^2 \leq
\tau \|\DeltaHat\|_1$, and hence
\begin{eqnarray*}
\|\DeltaHat\|_2^2 & \leq & \sqrt{2 \myrad} \tau^{1-\qpar/2}
\|\DeltaHat\|_2 + 2 \myrad \tau^{2-\qpar}.
\end{eqnarray*}
Viewed as a quadratic in the indeterminate $x = \|\DeltaHat\|_2$, this
inequality is equivalent to the constraint $f(x) = a x^2 + bx + c \leq
0$, with $a = 1$,
\begin{equation*}
b = - \sqrt{2 \myrad} \tau^{1 -\qpar/2}, \quad \mbox{and} \quad c = -2
\myrad \tau^{2-\qpar}.
\end{equation*}
Since $f(0) = c < 0$ and the positive root of $f(x)$ occurs at
$x^* = (-b + \sqrt{b^2 - 4ac})/(2a)$, some algebra shows that
we must have
\begin{equation*}
\|\DeltaHat\|_2^2 \; \leq \; 4 \max \{ b^2, \: |c| \} \; \leq \; 
24
\myrad \Big[ \frac{\colnormsq}{\conlowersq}
\frac{\sigma^2}{\conlowersq} \; \frac{\log \pdim}{\numobs}
\Big]^{1 - \qpar/2},
\end{equation*}
with high probability (stated in Theorem~\ref{ThmUpper}(a) which completes the proof of Theorem~\ref{ThmUpper}(a).

\myparagraph{Proof of Theorem~\ref{ThmUpper}(b)}
In order to establish the bound~\eqref{EqnUpperBoundZero}, we follow
the same steps with $\SpecFunLowerZ = 0$, thereby obtaining the following
simplified form of the bound~\eqref{EqnFirstQuad}:
\begin{eqnarray*}
\|\DeltaHat\|_2^2 & \leq & \frac{\colnorm}{\conlower} \;
\frac{\sigma}{\conlower} \; \sqrt{\frac{3 \log \pdim}{\numobs}} \,
\|\DeltaHat\|_1.
\end{eqnarray*}
By definition of the estimator, we have $\|\DeltaHat\|_0 \leq 2
\spindex$, from which we obtain $\|\DeltaHat\|_1 \leq \sqrt{2 \spindex}
\|\DeltaHat\|_2$.  Canceling out a factor of $\|\DeltaHat\|_2$ from
both sides yields the claim~\eqref{EqnUpperBoundZero}.

\vspace*{.1in}

Establishing the sharper upper bound~\eqref{EqnUpperBoundZeroSharp}
requires more precise control on the right-hand side of the
inequality~\eqref{EqnBasic}.  The following lemma, proved in
Appendix~\ref{AppLemNewrad}, provides this control:
%%%%%%%%%%%%%%%
\blems
\label{LemNewrad}
If $\frac{\|\Xmat \mytheta\|_2}{\sqrt{\numobs} \|\mytheta\|_2} \leq
\conupper$ for all $\mytheta \in \Ball_0(2 \s)$, then for any $\radtwo
> 0$, we have
\begin{eqnarray}
\label{EqnNewrad}
\sup_{\|\mytheta\|_0 \leq 2 \spindex, \|\theta\|_2 \leq \radtwo}
\frac{1}{\numobs} \big| \wnoise^T \Xmat \theta \big| & \leq & 6 \:
\sigma \; \radtwo \; \conupper \, \sqrt{\frac{\spindex
\log(\pdim/\spindex)}{\numobs}}
\end{eqnarray}
with probability greater than $1-\gencon_1 \exp(-\gencon_2 \min \{
\numobs, \spindex \log(\pdim - \spindex) \})$. 
\elems
Let us apply this lemma to the basic inequality~\eqref{EqnBasic}.  We may
upper bound the right-hand side as
\begin{eqnarray*}
\big| \frac{w^T \Xmat \Delta}{\numobs} \big| & \leq & \|\Delta\|_2 \;
\sup_{\|\mytheta\|_0 \leq 2 \spindex, \|\theta\|_2 \leq 1}
\frac{1}{\numobs} \big| \wnoise^T \Xmat \theta \big| \; \leq \; 6 \:
\|\Delta\|_2 \; \sigma \; \conupper \, \sqrt{\frac{\spindex
\log(\pdim/\spindex)}{\numobs}}.
\end{eqnarray*}
Consequently, we have
\begin{eqnarray*}
\frac{1}{\numobs} \| \Xmat \DeltaHat\|_2^2 & \leq & 12 \: \sigma \;
\|\DeltaHat\|_2 \; \conupper \, \sqrt{\frac{\spindex
\log(\pdim/\spindex)}{\numobs}},
\end{eqnarray*}
with high probability.  By Assumption~\ref{AssLower}, we have $\|\Xmat
\DeltaHat\|_2^2/\numobs \geq \conlowersq \|\DeltaHat\|^2_2$.
Cancelling out a factor of $\|\DeltaHat\|_2$ and re-arranging yields
$\|\DeltaHat\|_2 \leq 12 \: \frac{\conupper \sigma }{\conlowersq} \,
\sqrt{\frac{\spindex \log(\pdim/\spindex)}{\numobs}}$ with high probability as claimed.

\subsubsection{Proof of Theorem~\ref{ThmUpperPred}}

We again make use of the elementary inequality~\eqref{EqnBasic} to
establish upper bounds on the prediction risk.

\myparagraph{Proof of Theorem~\ref{ThmUpperPred}(a)}

So as to facilitate tracking of constants in this part of the proof,
we consider the rescaled observation model, in which $\wtil \sim N(0,
I_\numobs)$ and $\Xtil \defn
\sigma^{-1} \Xmat$. Note that if $\Xmat$ satisfies
Assumption~\ref{AssColNorm} with constant $\colnorm$, then $\Xtil$
satisfies it with constant $\colnormtil = \colnorm/\sigma$.  Moreover,
if we establish a bound on $\|\Xtil (\betahat -
\betastar)\|_2^2/\numobs$, then multiplying by $\sigma^2$ recovers a
bound on the original prediction loss.

We first deal with the case $\qpar = 1$.  In particular, we have
\begin{eqnarray*}
\big| \frac{1}{\numobs} \wtil^T \Xtil \mytheta \big| & \leq & \|
\frac{\wtil^T \Xtil}{\numobs} \|_\infty \|\mytheta\|_1 \; \leq \;
\sqrt{\frac{3 \colnormtilsq \sigma^2 \log \pdim}{\numobs}} \, (2 \, R_1),
\end{eqnarray*}
where the second inequality holds with probability $1-\gencon_1
\exp(-\gencon_2 \log \pdim)$, using standard Gaussian tail bounds.
(In particular, since $\|\Xtil_i\|_2/\sqrt{\numobs} \leq \colnormtil$,
the variate $\wtil^T \Xtil_i/\numobs$ is zero-mean Gaussian with
variance at most $\colnormtilsq/\numobs$.) This completes the proof
for $\qpar=1$.

Turning to the case $\qpar \in (0,1)$, in order to establish upper
bounds over $\Ballq(2 \myrad)$, we require the following analog of
Lemma~\ref{LemNewrad}, proved in Appendix~\ref{AppLemNewradPred}.  So
as to lighten notation, let us introduce the shorthand $\qRate \defn
\sqrt{\myrad} \; (\frac{\log
\pdim}{\numobs})^{\frac{1}{2}-\frac{\qpar}{4}}$.
\blems
\label{LemNewradPred}
For $\qpar \in (0,1)$, suppose that $\qRate = o(1)$ and $\pdim =
\Omega(\numobs)$.  Then for any fixed radius $\radtwo$ such that
$\radtwo \geq \gencon_3 \colnormtil^{\frac{\qpar}{2}} \, \qRate$ for
some numerical constant $\gencon_3 > 0$, we have
\begin{eqnarray*}
\sup_{\mytheta \in \Ballq(2 \myrad), \; \frac{\|\Xtil
\mytheta\|_2}{\sqrt{\numobs}} \leq \radtwo} \frac{1}{\numobs} \, \big|
\wtil^T \Xtil \theta \big| & \leq & \gencon_4 \radtwo \;
\colnormtil^{\frac{\qpar}{2}} \, \sqrt{\myrad} \; (\frac{\log
\pdim}{\numobs})^{\frac{1}{2} - \frac{\qpar}{4}},
\end{eqnarray*}
with probability greater than $1 - \gencon_1 \exp(-\gencon_2 \,
 \numobs \, \qRatesq)$.
\elems
Note that Lemma~\ref{LemNewradPred} above holds for any fixed radius
$\radtwo \geq \gencon_3 \colnormtil^{\frac{\qpar}{2}} \, \qRate$.  We
would like the apply the result of Lemma~\ref{LemNewradPred} to
$\radtwo = \frac{\|\Xmat \Delta\|_2}{\sqrt{n}}$, which is a random
quantity.  In Appendix~\ref{AppLemPeel}, we state and prove a
``peeling'' result that allows us to strengthen
Lemma~\ref{LemNewradPred} in a way suitable for our needs.  In
particular, if we define the event
\begin{eqnarray}
\label{EqnDefnKeyEvent}
\mathcal{E} & \defn & \big \{ \exists \; \theta \in \Ballq(2 \myrad)
\mbox{ such that } \frac{1}{\numobs} \, \big|\wtil^T \Xtil \theta
\big| \geq \gencon_4 \frac{\|\Xtil \theta \|_2}{\sqrt{n}}
\;\colnormtil^{\frac{\qpar}{2}} \, \sqrt{\myrad} \; (\frac{\log
\pdim}{\numobs})^{\frac{1}{2} - \frac{\qpar}{4}} \big \},
\end{eqnarray}
then we claim that
\begin{eqnarray*}
\mprob[\mathcal{E}] & \leq & \frac{2 \exp(-\gencon \, \numobs \,
 \qRatesq)}{1 - \exp(- \gencon \, \numobs \, \qRatesq)}.
\end{eqnarray*}
This claim follows from Lemma~\ref{LemPeel} in
Appendix~\ref{AppLemPeel} by making the choices $f_n(v; X_n) =
\frac{1}{n}| \w^T \Xmat v|$, $\rho(v) = \frac{\|\Xmat
v\|_2}{\sqrt{n}}$, and $g(\radtwo) = \gencon_3 \: \radtwo \:
\colnormtil^{\frac{\qpar}{2}} \, \sqrt{\myrad} \; (\frac{\log
\pdim}{\numobs})^{\frac{1}{2} - \frac{\qpar}{4}}$.

Returning to the main thread, from the basic
inequality~\eqref{EqnBasic}, when the event $\mathcal{E}$ from
equation~\eqref{EqnDefnKeyEvent} holds, we have
\begin{eqnarray*}
\frac{\|\Xtil \Delta\|_2^2}{\numobs} & \leq & % & \leq & \gencon \;
\frac{\|\Xtil \Delta\|_2}{\sqrt{\numobs}} \; \sqrt{\colnormtil^\qpar
\myrad \big(\frac{\log \pdim}{\numobs} \big)^{1- \qpar/2}}.
\end{eqnarray*}
Canceling out a factor of $\frac{\|\Xmat \Delta\|_2}{\sqrt{\numobs}}$,
squaring both sides, multiplying by $\sigma^2$ and simplifying yields
\begin{eqnarray*}
\frac{\|\Xmat \Delta\|_2^2}{\numobs} & \leq & \gencon^2 \; \sigma^2
 \big(\frac{\colnorm}{\sigma} \big)^\qpar \; \myrad \big(\frac{\log
 \pdim}{\numobs} \big)^{1- \qpar/2} \; = \; \gencon^2 \, \colnormsq \;
 \myrad \; \big(\frac{\sigma^2}{\colnormsq} \; \frac{\log
 \pdim}{\numobs} \big)^{1- \qpar/2},
\end{eqnarray*}
as claimed.

\myparagraph{Proof of Theorem~\ref{ThmUpperPred}(b)}

For this part, we require the following lemma, proven in
Appendix~\ref{AppLemNewradPredZero}:
\blems
\label{LemNewradPredZero}
Suppose that $\frac{\pdim}{2 \spindex} \geq 2$.  Then for any $\radtwo
> 0$, we have
\begin{eqnarray*}
\sup_{\mytheta \in \Ball_0(2 \spindex), \frac{\|\Xmat
\mytheta\|_2}{\sqrt{\numobs}} \leq \radtwo} \frac{1}{\numobs} \big|
\wnoise^T \Xmat \mytheta \big| & \leq & 9 \; \radtwo \: \sigma \;
\sqrt{\frac{\spindex \log(\frac{\pdim}{\spindex})}{\numobs}}
\end{eqnarray*}
with probability greater than $1- \exp \big( -10 \spindex
\log(\frac{\pdim}{2 \spindex}) \big)$.
\elems
Consequently, combining this result with the basic
inequality~\eqref{EqnBasic}, we conclude that 
\begin{align*}
\frac{\|\Xmat \Delta\|^2_2}{\numobs} & \leq 9 \, \frac{\|\Xmat
\Delta\|_2}{\sqrt{\numobs}} \, \sigma \; \sqrt{\frac{\spindex
\log(\frac{\pdim}{\spindex})}{\numobs}},
\end{align*}
with high probability, from which the result follows.

\section{Discussion}
\label{Discussion}

The main contribution of this paper was to analyze minimax rates of
convergence for the linear model~\eqref{EqnLinearObs} under
high-dimensional scaling, in which the sample size $\numobs$ and
problem dimension $\pdim$ tend to infinity. We provided lower bounds
for the $\ell_{\rpar}$-norm for all $\rpar \in [1, \infty]$ with
$\rpar \neq \qpar$, as well as for the $\LTwo$-prediction loss.  In
addition, for both the $\ell_2$-loss and $\Ltwo$-prediction loss, we
derived a set of upper bounds that match our lower bounds up to
constant factors, so that the minimax rates are exactly determined in
these cases. The rates may be viewed as an extension of the rates for
the case of $\ell_2$-loss from Donoho and Johnstone~\cite{DonJoh94} on
the Gaussian sequence model to more general design matrices
$\Xmat$. In particular substituting $\Xmat=I$ and $\pdim=n$ into
Theorems~\ref{ThmLower} and \ref{ThmUpper}, yields the same rates as
those expressed in Donoho and Johnstone~\cite{DonJoh94} (see
Corollary~\ref{CorNorSeqModel}), although they provided much sharper
control of the constant pre-factors than the analysis given here.

Apart from the rates themselves, our analysis highlights how
conditions on the design matrix $\Xmat$ enter in complementary manners
for different loss functions.  On one hand, it is possible to obtain
lower bounds on $\ell_2$-risk (see Theorem~\ref{ThmLower}) or upper
bounds on $\Ltwo$-prediction risk (see Theorem~\ref{ThmUpperPred})
under very mild assumptions on $\Xmat$---in particular, our analysis
requires only that the columns of $\Xmat/\sqrt{\numobs}$ have bounded
$\ell_2$-norms (see, in particular, Assumption~\ref{AssColNorm}).  On
the other hand, in order to obtain upper bounds on $\ell_2$ risk
(Theorem~\ref{ThmUpper}) or lower bound on $\LTwo$-norm prediction
risk (Theorem~\ref{ThmLowerPred}), the design matrix $\Xmat$ must
satisfy, in addition to column normalization, other more restrictive
conditions.  In particular, our analysis was based on imposed on a
certain type of lower bound on the curvature of $\Xmat^T \Xmat$
measured over the $\ell_\qpar$-ball (see Assumption~\ref{AssLower}).
As shown in Lemma~\ref{LemIdent}, this lower bound is intimately
related to the \emph{degree of non-identifiability} over the $\ell_{\qpar}$-ball of the high-dimensional linear regression model .

In addition, we showed that Assumption~\ref{AssLower} is not
unreasonable---in particular, it is satisfied with high probability
for broad classes of Gaussian random matrices, in which each row is
drawn in an i.i.d. manner from a $N(0, \CovMat)$ distribution (see
Proposition~\ref{PropRandDesign}).  This result applies to Gaussian
ensembles with much richer structure than the standard Gaussian case
($\CovMat = I_{\pdim \times \pdim}$).  Finally, we compared to the
weakest known sufficient conditions for $\ell_1$-based relaxations to
be consistent in $\ell_2$-norm for $\qpar=0$---namely, the restricted
eigenvalue (RE) condition, of Bickel et al.~\cite{BicRitTsy08} and
showed that the oracle least-squares over the $\ell_0$-ball method can
succeed with even milder conditions on the design. In addition, we
also proved that the RE condition holds with high probability for
broad classes for Gaussian random matrices, as long as the covariance
matrix $\CovMat$ is not degenerate. The analysis highlights how the
structure of $\Xmat$ determines whether $\ell_1$-based relaxations
achieve the minimax optimal rate.

The results and analysis from our paper can be extended in a number of
ways. First, the assumption of independent Gaussian noise is somewhat
restrictive and it would be interesting to analyze the model under
different noise assumption, either noise with heavier tails or some
degree of dependency. In addition, we are currently working on
extending our analysis to non-parametric sparse additive models.

\subsection*{Acknowledgements}

We thank Arash Amini for useful discussion, particularly regarding the
proofs of Proposition~\ref{PropRandDesign} and Lemma~\ref{LemPeel}.
This work was partially supported by NSF grant DMS-0605165 to MJW and
BY.  In addition, BY was partially supported by the NSF grant
SES-0835531 (CDI), the NSFC grant 60628102 and a grant from the MSRA.
MJW was supported by an Sloan Foundation Fellowship and AFOSR Grant
FA9550-09-1-0466.  During this work, GR was financially supported by a
Berkeley Graduate Fellowship.

%%%%%%%%%%%%%%%%%%%%%%%%%%%%%%%%%%%%%%%%%%%%%%%%%%%%%%%%%%%%%%%%%%%%%%%%%

\appendix

%%%%%%%%%%%%%%%%%%%%%%%%%%%%%%%%%%%%%%%%%%%%%%%%%%%%%%%%%%%%%%%%%%%%%%%%%%

\section{Proof of Lemma 2}
\label{AppLemElloneElltwo}
Defining the set $S = \{ j \, \mid |\mytheta_j| > \tau \}$, we have
\begin{eqnarray*}
\|\mytheta \|_1 & = & \|\mytheta_S\|_1 + \sum_{j \notin S}
|\mytheta_j| \; \leq \; \sqrt{|S|} \|\mytheta\|_2 + \tau \sum_{j
\notin S} \frac{|\mytheta_j|}{\tau}.
\end{eqnarray*}
Since $|\mytheta_j|/\tau < 1$ for all $i \notin S$, we obtain
\begin{eqnarray*}
\|\mytheta\|_1 & \leq & \sqrt{|S|} \|\mytheta\|_2 + \tau \sum_{j
\notin S} \big(|\mytheta_i|/\tau\big)^\qpar \\
& \leq & \sqrt{|S|} \|\mytheta\|_2 + 2 \myrad \tau^{1-\qpar}.
\end{eqnarray*}
Finally, we observe $2 \myrad \geq \sum_{j \in S} |\mytheta_j|^\qpar
\; \geq \; |S| \tau^\qpar$, from which the result follows.

\section{Proof of Lemma~\REFLEMENTROPY}
\label{AppLemEntropyBounds}

The result is obtained by inverting known results on (dyadic) entropy
numbers of $\ell_\qpar$-balls; there are some minor technical
subtleties in performing the inversion.  For a $\pdim$-dimensional
$\ell_{\qpar}$ ball with $\qpar \in (0,\rpar)$, it is
known~\cite{Sch84,Kuhn,GueLit00} that for all integers \mbox{$k \in
[\log \pdim, \pdim]$,} the dyadic entropy numbers $\dyaeps_k$ of the
ball $\Ballq(1)$ with respect to the $\ell_\rpar$-norm scale as
\begin{equation}
\label{EqnSchutt}
\dyaeps_k(\ell_{\qpar}, \|\cdot\|_\rpar) = \Cqr \, \biggr[ \frac{\log
(1 + \frac{\pdim}{k})}{k} \biggr]^{1/\qpar - 1/\rpar}.
\end{equation}
Moreover, for $k \in [1, \log \pdim]$, we have $\dyaeps_k(\ell_\qpar)
\leq \Cqr$.

We first establish the upper bound on the metric entropy.  Since
$\pdim \geq 2$, we have
\begin{eqnarray*}
\entnum_k(\ell_q) & \leq & \Cqr \biggr[ \frac{\log
(1+\frac{\pdim}{2})}{ k} \biggr]^{1/\qpar - 1/\rpar} \; \leq \; \;
\Cqr \biggr[ \frac{\log \pdim}{ k} \biggr]^{1/\qpar - 1/\rpar}.
\end{eqnarray*}
Inverting this inequality for $k = \log \CovNum_{\rpar,
\qpar}(\epsilon)$ and allowing for a ball radius $\myrad$ yields
\begin{eqnarray}
\label{EqnUpperLemma}
\log \CovNum_{\rpar, \qpar}(\epsilon) & \leq & \big( \Cqr \, \frac
{\Rq^{1/\qpar}}{\epsilon} \big)^{\frac{\rpar \qpar}{\rpar-\qpar}} \log
\pdim,
\end{eqnarray}
as claimed.

We now turn to proving the lower bound on the metric entropy, for
which we require the existence of some fixed $\kOne \in (0,1)$ such
that $\kdim \leq \pdim^{1-\kOne}$. Under this assumption, we have $1 +
\frac{\pdim}{\kdim} \geq \frac{\pdim}{\kdim} \geq \pdim^\kOne$, and
hence
\begin{eqnarray*}
\Cqr \biggr[ \frac{\log(1 + \frac{\pdim}{k})}{k} \biggr]^{1/\qpar -
1/\rpar} & \geq & \Cqr \biggr[\frac{\kOne \log \pdim
}{k}\biggr]^{1/\qpar - 1/\rpar}\\
\end{eqnarray*}
Accounting for the radius $\Rq$ as was done for the upper bound yields
\begin{eqnarray*}
\log \CovNum_{\rpar, \qpar}(\epsilon) & \geq & \kOne \big(\frac{\Cqr
\Rq^{1/\qpar}}{\epsilon} \big)^{\frac{\rpar \qpar}{\rpar-\qpar}} \log
\pdim,
\end{eqnarray*}
as claimed.

Finally, let us check that our assumptions on $\kdim$ needed to
perform the inversion are ensured by the conditions that we have
imposed on $\epsilon$. The condition $\kdim \geq \log \pdim$ is
ensured by setting $\epsilon < 1$. Turning to the condition $\kdim
\leq \pdim^{1-\kOne}$, from the bound~\eqref{EqnUpperLemma} on
$\kdim$, it suffices to choose $\epsilon$ such that $\big(
\frac{\Cqr}{\epsilon} \big)^{\frac{\rpar \qpar}{\rpar-\qpar}} \: \log
\pdim \leq \pdim^{1-\kOne}$. This condition is ensured by enforcing
the lower bound $\epsilon^\rpar = \Omega\big(\frac{\log
\pdim}{\pdim^{1-\kOne}} \big)^{\frac{\rpar- \qpar}{\qpar}}$ for some
$\kOne \in (0,1)$.

% MY VOTE IS TO IGNORE THE EXPLICIT CONSTANT; IT MAKES THINGS SEEM
% MORE COMPLICATED WITHOUT ADDING MUCH (MJW)
%\Gue and Litvak~\cite{GueLit00} give an explicit upper bound on the
%constant $\Cqr$, which yields
%\begin{eqnarray*}
%\log \CovNum_{\rpar, \qpar}(\epsilon) & \leq & & \leq & \UqrOne
% \big(\frac{\Rq^{1/q}}{\epsilon} \big)^{\frac{\rpar
% \qpar}{\rpar-\qpar}} \log \pdim \nonumber
%\end{eqnarray*}
%where
%\begin{eqnarray*}
%\UqrOne & = & 2^{\frac{\rpar+1}{\rpar-\qpar}}
%\biggr(\frac{c_1}{q'}\biggr)^{\frac{\rpar \qpar}{\rpar-\qpar}}
%\log\big(\frac{2}{q'}\big).
%\end{eqnarray*}

%%%%%%%%%%%%%%%%%%%%%%%%%%%%%%%%%%%%%%%%%%%%%%%%%%%%%%%%%%%%%%%%%%%%%%%%

\section{Proof of Lemma~\REFLEMKL}
\label{AppLemKL}
We deal first with (dyadic) entropy numbers, as previously
defined~\eqref{EqnDefnDyaEps}, and show that
\begin{eqnarray}
\label{EqnInterDya} 
\dyaeps_{2k-1}(\absconv_\qpar(\Xmat/\sqrt{\numobs}), \: \|\cdot\|_2) &
\leq & \gencon \, \colnorm \; \min \biggr \{1, \big(\frac{\log(1 +
\frac{\pdim}{k})}{k} \big)^{\frac{1}{\qpar} - \frac{1}{2}} \biggr \}.
\end{eqnarray}
We prove this intermediate claim by combining a number of known
results on the behavior of dyadic entropy numbers.  First, using
Corollary 9 from \Gue and Litvak~\cite{GueLit00}, for all $k = 1, 2,
\ldots$, we have
\begin{eqnarray*}
\dyaeps_{2k-1}(\absconv_\qpar(\Xmat/\sqrt{\numobs}), \: \|\cdot\|_2) &
\leq & \gencon \; \dyaeps_k(\absconv_1(\Xmat), \: \|\cdot\|_2) \; \min
\biggr \{1, \big(\frac{\log(1 + \frac{\pdim}{k})}{k}
\big)^{\frac{1}{\qpar}-1} \biggr \}.
\end{eqnarray*}
Using Corollary 2.4 from Carl and Pajor~\cite{CarPaj88}, we obtain
\begin{eqnarray*}
 \dyaeps_k(\absconv_1(\Xmat/\sqrt{\numobs}), \: \|\cdot\|_2) & \leq &
\frac{\gencon}{\sqrt{\numobs}} \; \matsnorm{\Xmat}{1 \rightarrow 2} \;
\min \biggr \{1, \big(\frac{\log(1 + \frac{\pdim}{k})}{k} \big)^{1/2}
\biggr \},
\end{eqnarray*}
where $\matsnorm{\Xmat}{1 \rightarrow 2}$ denotes the norm of $\Xmat$
viewed as an operator from $\ell_1^\pdim \rightarrow \ell_2^\numobs$.
More specifically, we have
\begin{eqnarray*}
\frac{1}{\sqrt{\numobs}} \matsnorm{\Xmat}{1 \rightarrow 2} & = &
\frac{1}{\sqrt{\numobs}} \, \sup_{\|u\|_1 = 1} \| \Xmat u\|_2 \\
& = & \frac{1}{\sqrt{\numobs}} \sup_{\|v\|_2 = 1} \sup_{\|u\|_1 = 1}
v^T X u \\
& = & \max_{i = 1, \ldots, \pdim} \|\Xmat_i\|_2/\sqrt{\numobs} \; \leq
\; \colnorm.
\end{eqnarray*}
Overall, we have shown that
$\dyaeps_{2k-1}(\absconv_\qpar(\Xmat/\sqrt{\numobs}), \: \|\cdot\|_2)
\, \leq \, \gencon \, \colnorm \; \min \biggr \{1, \big(\frac{\log(1 +
\frac{\pdim}{k})}{k} \big)^{\frac{1}{\qpar} - \frac{1}{2}} \biggr \}$,
as claimed.  Finally, under the stated assumptions, we may invert the
upper bound~\eqref{EqnInterDya} by the same procedure as in the proof
of Lemma~\ref{LemEntropyBounds} (see
Appendix~\ref{AppLemEntropyBounds}), thereby obtaining the claim.

%%%%%%%%%%%%%%%%%%%%%%%%%%%%%%%%%%%%%%%%%%%%%%%%%%%%%%%%%%%%%%%%%%%%

\section{Proof of Lemma~\REFLEMHYPER}
\label{AppLemHyperCube}

In this appendix, we prove Lemma~\ref{LemHyperCube}.  Our proof is
inspired by related results from the approximation theory literature
(see, e.g., K\"{u}hn~\cite{Kuhn}).  For each even integer $\spindex =
2, 4, 6, \ldots, \pdim$, let us define the set
\begin{equation}
\Hyper \defn \big \{ z \in \{-1, 0, +1\}^\pdim \, \mid \, \|z\|_0 =
\spindex \big \}.
\end{equation}
Note that the cardinality of this set is $|\Hyper| = {\pdim \choose
\spindex} 2^{\spindex}$, and moreover, we have $\|z - z'\|_0 \leq 2
\spindex$ for all pairs $z, z' \in \Hyper$.  We now define the Hamming
distance $\rhohamm$ on $\Hyper \times \Hyper$ via \mbox{$\rhohamm(z,
z') = \sum_{j=1}^\pdim \Ind[z_j \neq z'_j]$.} For some fixed element
$z \in \Hyper$, consider the set $\{ z' \in \Hyper \, \mid \,
\rhohamm(z, z') \leq s/2 \}$.  Note that its cardinality is upper
bounded as
\begin{equation*}
\big | \{ z' \in \Hyper \, \mid \, \rhohamm(z, z') \leq s/2 \} \big|
\; \leq \; {\pdim \choose \s/2 } 3^{\s/2}.
\end{equation*}
To see this, note that we simply choose a subset of size $\s/2$ where
$z$ and $z'$ agree and then choose the other $\s/2$ co-ordinates
arbitrarily.

Now consider a set $\Aset \subset \Hyper$ with cardinality at most
$|\Aset| \leq m \defn \frac{{\pdim \choose \s}}{{\pdim \choose
\s/2}}$.  The set of elements $z \in \Hyper$ that are within Hamming
distance $\spindex/2$ of some element of $\Aset$ has cardinality at
most
\begin{eqnarray*} |\{z \in \Hyper \, \mid \, : \rhohamm(z,z')\leq
\spindex/2 \mbox{ for some $z' \in \Aset$} \}| & \leq & |\Aset| \;
{\pdim \choose \spindex/2} 3^{\s/2 } \; < \; |\Hyper|,
\end{eqnarray*}
where the final inequality holds since $ m {\pdim \choose \spindex/2}
3^{\s/2 } \; < \; |\Hyper|$.  Consequently, for any such set with
cardinality $|\Aset| \leq m$, there exists a $z \in \Hyper$ such that
$\rhohamm(z, z') > \spindex/2$ for all $z' \in \Aset$.  By inductively
adding this element at each round, we then create a set with $\Aset
\subset \Hyper$ with $|\Aset| > m$ such that $\rhohamm(z, z') >
\spindex/2$ for all $z, z' \in \Aset$.

To conclude, let us lower bound the cardinality $m$.  We have
\begin{equation*}
m = \frac{{\pdim \choose \s}}{{\pdim \choose \s/2}} \; = \;
\frac{(\pdim - \spindex/2)! \, (\spindex/2)!}{(\pdim - \spindex)! \;
\spindex!}  = \prod_{j=1}^{\s/2} \frac{\pdim -\s+j}{\s/2+j} \; \geq \;
\big(\frac{\pdim - \spindex}{\spindex/2} \big)^{\spindex/2},
\end{equation*}
where the final inequality uses the fact that the ratio $\frac{\pdim
  -\s+j}{\s/2+j}$ is decreasing as a function of $j$.

%%%%%%%%%%%%%%%%%%%%%%%%%%%%%%%%%%%%%%%%%%%%%%%%%%%%%%%%%%%%%%%%%%%%%%%%%%

%%%%%%%%%%%%%%%%%%%%%%%%%%%%%%%%%%%%%%%%%%%%%%%%%%%%%%%%%%%%%%%%%%%%%%%%%%

\section{Proof of Proposition 1}
\label{AppPropRandDesign}

In this appendix, we prove both parts of
Proposition~\ref{PropRandDesign}.  In addition to proving the lower
bound~\eqref{EqnLowerFinal}, we also prove the analogous upper bound
\begin{align}
\label{EqnUpperFinal}
\frac{\|X v\|_2}{\sqrt{\numobs}} & \leq 3 \|\CovMatSqrt v\|_2 + 6
\biggr[\frac{\rhomax(\CovMat)\log \pdim}{\numobs} \biggr]^{1/2}
\|v\|_1 \qquad \mbox{ for all $v \in \real^\pdim$.}
\end{align}
Our approach to proving the bounds~\eqref{EqnLowerFinal}
and~\eqref{EqnUpperFinal} is based on Slepian's
lemma~\cite{LedTal91,DavSza01} as well as an extension thereof due to
Gordon~\cite{Gordon85}.  For the reader's convenience, we re-state
versions of this lemma here.  Given some index set $U \times V$, let
$\{Y_{u,v}, \: (u,v) \in U \times V \}$ and $\{Z_{u,v}, \: (u,v) \in U
\times V\}$ be a pair of zero-mean Gaussian processes.  Given the
semi-norm on these processes defined via $\sigma(X) =
\Exs[X^2]^{1/2}$, Slepian's lemma asserts that if
\begin{eqnarray}
\label{EqnSlepianCondition}
\sigma(Y_{u,v}-Y_{u',v'}) & \leq & \sigma(Z_{u,v} - Z_{u',v'}) \qquad
\mbox{ for all $(u,v)$ and $(u',v')$ in $U \times V$,}
\end{eqnarray}
then 
\begin{align}
\label{EqnSlepian}
\Exs[\sup_{(u,v) \in U \times V} Y_{u,v}] & \leq \Exs[\sup_{(u,v) \in
U \times V} Z_{u,v}].
\end{align}
One version of Gordon's extension~\cite{Gordon85,LedTal91} asserts
that if the inequality~\eqref{EqnSlepianCondition} holds for $(u,v)$
and $(u', v')$ in $U \times V$, and holds with \emph{equality} when $v
= v'$, then
\begin{eqnarray}
\label{EqnGordon}
\Exs[\sup_{u \in U}\; \inf_{v \in V} Y_{u,v}] & \leq & \Exs[\sup_{u \in
U}\; \inf_{v \in V} Z_{u,v}].
\end{eqnarray}

\newcommand{\SuperSpecSet}[1]{\ensuremath{V_{}(\newrad)}}
\newcommand{\Sphere}[1]{\ensuremath{S^{#1-1}}}

Turning to the problem at hand, any random matrix $X$ from the given
ensemble can be written as $W \CovMatSqrt$, where $W \in
\real^{\numobs \times \pdim}$ is a matrix with i.i.d. $N(0,1)$
entries, and $\CovMatSqrt$ is the symmetric matrix square root.  We
choose the set $U$ as the unit ball $\Sphere{\numobs} = \{ u \in
\real^{\numobs} \, \mid \, \|u\|_2 = 1 \}$, and for some radius
$\newrad$, we choose $V$ as the set
\begin{align*}
\SuperSpecSet{1} & \defn  \{ v \in \real^\pdim \, \mid \, \|\CovMatSqrt
v\|_2 = 1, \|v\|_\qpar^\qpar \leq \newrad \}.
\end{align*}
(Although this set may be empty for certain choices of $\newrad$, our
analysis only concerns those choices for which it is non-empty.)  For
a matrix $M$, we define the associated Frobenius norm
\mbox{$\matsnorm{M}{F} = [\sum_{i,j} M_{ij}^2]^{1/2}$,} and for any $v
\in \SuperSpecSet{\qpar}$, we introduce the convenient shorthand
$\vtil = \CovMatSqrt \: v$.

With these definition, consider the centered Gaussian process $Y_{u,v}
= u^T W v$ indexed by $\Sphere{\numobs} \times \SuperSpecSet{1}$.
Given two pairs $(u,v)$ and $(u', v')$ in $\Sphere{\numobs} \times
\SuperSpecSet{1}$, we have
\begin{eqnarray}
\sigma^2(Y_{u,v} - Y_{u',v'}) & = & \matsnorm{u \, \vtil^T - u'
(\vtil')^T}{F}^2 \nonumber \\
& = & \matsnorm{u \vtil^T - u' \, \vtil^T + u'\, \vtil^T - u'
(\vtil')^T}{F}^2 \nonumber \\
\label{EqnKeyDom}
& = & \|\vtil\|_2^2 \, \|u-u'\|_2^2 + \|u'\|_2^2 \|\vtil -
\vtil'\|_2^2 + 2 (u^T u' - \|u'\|_2^2) (\|\vtil\|_2^2 - \vtil^T
\vtil') 
\end{eqnarray}
Now by the Cauchy-Schwarz inequality and the equalities $\|u\|_2 =
\|u'\|_2 = 1$ and $\|\vtil\|_2 = \|\vtil'\|_2$, we have \mbox{$u^T u'
- \|u\|_2^2 \leq 0$,} and $\|\vtil\|_2^2 - \vtil^T \vtil' \geq 0$.
Consequently, we may conclude that
\begin{eqnarray}
\label{EqnKeyDomB}
\sigma^2(Y_{u,v} - Y_{u',v'}) & \leq & \|u-u'\|_2^2 + \|\vtil -
\vtil'\|_2^2.
\end{eqnarray}
We claim that the Gaussian process $Y_{u,v}$ satisfies the conditions
Gordon's lemma in terms of the zero-mean Gaussian process $Z_{u,v}$
given by
\begin{eqnarray}
Z_{u,v} & = &  g^T u + h^T
(\CovMatSqrt \, v),
\end{eqnarray}
where $g \in \real^n$ and $h \in \real^\pdim$ are both standard
Gaussian vectors (i.e., with i.i.d. $N(0,1)$ entries). To establish
this claim, we compute
\begin{eqnarray*}
\sigma^2(Z_{u,v} - Z_{u',v'}) & = & \|u - u'\|_2^2 + \|\CovMatSqrt \,
(v - v')\|_2^2 \\
& = & \|u - u'\|_2^2 + \|\vtil - \vtil'\|_2^2.
\end{eqnarray*}
Thus, from equation~\eqref{EqnKeyDomB}, we see that Slepian's
condition~\eqref{EqnSlepianCondition} holds.  On the other hand, when
$v = v'$, we see from equation~\eqref{EqnKeyDom} that
\begin{equation*}
\sigma^2(Y_{u,v} - Y_{u',v}) \, = \, \|u - u'\|_2^2 \, = \,
\sigma^2(Z_{u,v} - Z_{u,v'}),
\end{equation*}
so that the equality required for Gordon's inequality is also
satisfied.

\paragraph{Establishing an upper bound:}  We begin by
exploiting Slepian's inequality~\eqref{EqnSlepian} to establish the
upper bound~\eqref{EqnUpperFinal}.  We have
\begin{eqnarray*}
\Exs \big[ \sup_{v \in \SuperSpecSet{\qpar}} \|X v\|_2 \big] & = &
\Exs \big[ \sup_{(u,v) \in \Sphere{\numobs} \times
\SuperSpecSet{\qpar}} u^T X v] \\
& \leq & \Exs \big[ \sup_{(u,v) \in \Sphere{\numobs} \times
\SuperSpecSet{\qpar}} Z_{u,v} \big] \\
& = & \Exs [\sup_{\|u\|_2 = 1} g^T u] + \Exs[\sup_{v \in
\SuperSpecSet{\qpar}} h^T (\CovMatSqrt v)] \\
& \leq & \Exs [\|g\|_2] + \Exs[\sup_{v \in \SuperSpecSet{\qpar}} h^T
(\CovMatSqrt v)].
%& \leq & \Exs [\|g\|_2] + \newrad \; \Exs[\|\CovMatSqrt h\|_\infty],
\end{eqnarray*}
By convexity, we have $\Exs[\|g\|_2] \leq \sqrt{\Exs[\|g\|_2^2]} =
\sqrt{\numobs}$, from which we can conclude that
\begin{eqnarray}
\label{EqnEarly}
\Exs \big[ \sup_{v \in \SuperSpecSet{\qpar}} \|X v\|_2 \big] & \leq &
\sqrt{\numobs} + \Exs[\sup_{v \in \SuperSpecSet{\qpar}} h^T
(\CovMatSqrt v)].
\end{eqnarray}
Turning to the remaining term, we have
\begin{equation*}
\sup_{v \in \SuperSpecSet{1}} |h^T (\CovMatSqrt v)| \leq \sup_{v \in
\SuperSpecSet{1}} \|v\|_1 \; \| \CovMatSqrt h\|_\infty \; \leq \;
\newrad \| \CovMatSqrt h\|_\infty.
\end{equation*}
Since each element $(\CovMatSqrt h)_i$ is zero-mean Gaussian with
variance at most $\rhomax(\CovMat) = \max_{i} \CovMat_{ii}$, standard
results on Gaussian maxima (e.g.,~\cite{LedTal91}) imply that
$\Exs[\|\CovMatSqrt h\|_\infty] \leq \sqrt{3 \rhomax(\CovMat) \, \log
\pdim}$.  Putting together the pieces, we conclude that for $\qpar =
1$
\begin{eqnarray}
\label{EqnUpperExp}
\Exs \big[ \sup_{v \in \SuperSpecSet{1}} \|X v\|_2/\sqrt{\numobs}
\big] & \leq & \underbrace{1 + \big [3 \rhomax(\CovMat) \, \frac{\log
\pdim}{\numobs} \big]^{1/2} \; \newrad}. \\
& & \qquad \qquad \tupper(\newrad) \nonumber
\end{eqnarray}

Having controlled the expectation, it remains to establish sharp
concentration.  Let \mbox{$f: \real^D \rightarrow \real$} be Lipschitz
function with constant $L$ with respect to the $\ell_2$-norm.  Then if
\mbox{$w \sim N(0, I_{D \times D})$} is standard normal, we are
guaranteed~\cite{Ledoux01} that for all $t > 0$,
\begin{eqnarray}
\label{EqnLedoux}
\mprob \big[ |f(w) - \Exs[f(w)]| \geq t] & \leq & 2 \exp(-\frac{t^2}{2
L^2}).
\end{eqnarray}
Note the dimension-independent nature of this inequality.  We apply
this result to the random matrix $W \in \real^{\numobs \times \pdim}$,
viewed as a standard normal random vector in $D = \numobs \, \pdim$
dimensions.  First, letting $f(W) = \sup_{v \in \SuperSpecSet{\qpar}}
\|W \CovMatSqrt v\|_2/\sqrt{\numobs}$, we find that
\begin{eqnarray*}
\sqrt{\numobs} \big[f(W) - f(W') \big] & = & \sup_{v \in
 \SuperSpecSet{\qpar}} \|W \CovMatSqrt v\|_2 - \sup_{v \in
 \SuperSpecSet{\qpar}} \|W' \CovMatSqrt v\|_2 \\
& \leq & \sup_{v \in \SuperSpecSet{\qpar}} \|\CovMatSqrt v\|_2
\matsnorm{(W-W')}{F} \\
& = & \matsnorm{W - W'}{F}
\end{eqnarray*}
since $\|\CovMatSqrt v\|_2 = 1$ for all $v \in \SuperSpecSet{1}$.  We
have thus shown that the Lipschitz constant $L \leq 1/\sqrt{\numobs}$.
Recalling the definition of $\tupper(\newrad)$ from the upper
bound~\eqref{EqnUpperExp}, we set $t = \tupper(\newrad)/2$ in the tail
bound~\eqref{EqnLedoux}, thereby obtaining
\begin{eqnarray}
\label{EqnUpperConc}
\mprob \big[ \sup_{v \in \SuperSpecSet{\qpar}} \|X v\|_2 \geq
\frac{3}{2} \tupper(\newrad; \qpar) \big] & \leq & 2 \exp(-\numobs
\frac{\tupper(\newrad)^2}{8}).
\end{eqnarray}
We now exploit this family of tail bounds to upper bound the 
probability of the event
\begin{eqnarray*}
\Tail & \defn & \big \{ \exists \; v \in \real^\pdim \;
\mbox{s.t. $\|\CovMatSqrt v\|_2 = 1$ and } \|X v \|_2 \geq 3
\tupper(\|v\|_1) \big \}.
\end{eqnarray*}
We do so using Lemma~\ref{LemPeel} from Appendix~\ref{AppLemPeel}.  In
particular, for the case $\Event = \Tail$, we may apply this lemma
with the objective functions $f(v; X) = \|X v\|_2$, sequence
$a_\numobs = \numobs$, the constraint $\rho(\cdot) = \|\cdot\|_1$, the
set $S = \{v \in \real^\pdim \, \mid \, \|\CovMatSqrt v\|_2 = 1\}$,
and $g(\newrad) = 3 \tupper(\newrad)/2$.  Note that the
bound~\eqref{EqnUpperConc} means that the tail
bound~\eqref{EqnBasicTail} holds with $c = 4/72$.  Therefore, by
applying Lemma~\ref{LemPeel}, we conclude that $\mprob[\Tail] \leq
\plaincon_1 \exp(- \plaincon_2 \numobs)$ for some numerical constants
$\plaincon_i$.

Finally, in order to extend the inequality to arbitrary $v \in
\real^\pdim$, we note that the rescaled vector $\vbre =
v/\|\CovMatSqrt v\|_2$ satisfies $\|\CovMatSqrt \vbre\|_2 = 1$.
Consequently, conditional on the event $\Tail^c$, we have
\begin{eqnarray*}
\|X \vbre \|_2/\sqrt{\numobs} & \leq & 3 + 3 \big[\sqrt{(3
\rhomax(\CovMat) \,\log \pdim)/\numobs} \big] \; \|\vbre\|_1,
\end{eqnarray*}
or equivalently, after multiplying through by $\|\CovMatSqrt v\|_2$,
the inequality
\begin{eqnarray*}
\|X v\|_2/\sqrt{\numobs} & \leq & 3 \|\CovMatSqrt v\|_2 + 3 (\sqrt{(3
\rhomax(\CovMat) \,\log \pdim)/\numobs} ) \|v\|_1,
\end{eqnarray*}
thereby establishing the claim~\eqref{EqnUpperFinal}.

%%%%%%%%%%%%%%%%%%%%%%%%%%%%%%%%%%%%%%%%%%%%%%%%%%%%%%%%%%%%%%%%%%%%%%%%%

\paragraph{Establishing the lower bound~\eqref{EqnLowerFinal}:}  

We now exploit Gordon's inequality in order to establish the lower
bound~\eqref{EqnLowerFinal}.  We have
\begin{eqnarray*}
- \inf_{v \in \SuperSpecSet{\qpar}} \|X v\|_2 & = & \sup_{v \in V} -\|
X v\|_2 \; = \; \sup_{v \in \SuperSpecSet{\qpar}} \inf_{u \in U} u^T X
v.
\end{eqnarray*}
Applying Gordon's inequality, we obtain
\begin{eqnarray*}
\Exs[ \sup_{v \in \SuperSpecSet{\qpar}} -\|X v\|_2] & \leq & \Exs
    \big[\sup_{v \in \SuperSpecSet{\qpar}} \inf_{u \in
    \Sphere{\numobs}} Z_{u,v} \big] \\
  & = & \Exs[\inf_{u \in \Sphere{\numobs}} g^T u] + \Exs[\sup_{v \in
  \SuperSpecSet{\qpar}} h^T \CovMatSqrt v] \\
& \leq & - \Exs[\|g\|_2] + \big [3 \rhomax(\CovMat) \, \log
\pdim \big]^{1/2} \; \newrad.
\end{eqnarray*}
where we have used our previous derivation to upper bound
$\Exs[\sup_{v \in \SuperSpecSet{\qpar}} h^T \CovMatSqrt v]$.
Noting\footnote{In fact, $|\Exs[\|g\|_2] - \sqrt{\numobs}| =
o(\sqrt{\numobs})$, but this simple bound is sufficient for our
purposes.} that $\Exs[\|g\|_2] \geq \sqrt{\numobs}/2$ for all $\numobs
\geq 1$, we divide by $\sqrt{\numobs}$ and add $1$ to both sides so as
to obtain
\begin{eqnarray}
\label{EqnLowerExpAlt}
\Exs \big[ \sup_{v \in \SuperSpecSet{\qpar}} \big(1 -\|X
v\|_2/\sqrt{\numobs}\big) \big] & \leq & \underbrace{1/2 + 
\big [3 \rhomax(\CovMat) \, \log
\pdim \big]^{1/2} \; \newrad} \\
& & \qquad \qquad \tlower(\newrad) \nonumber
\end{eqnarray}

Next define the function $f(W) = \sup_{v \in \SuperSpecSet{\qpar}}
\big(1 -\|W \CovMatSqrt v\|_2/\sqrt{\numobs} \big)$.  The same
argument as before shows that its Lipschitz constant is at most
$1/\sqrt{\numobs}$.  Setting $t = \tlower(\newrad)/2$ in the
concentration statement~\eqref{EqnLedoux} and combining with the lower
bound~\eqref{EqnLowerExpAlt}, we conclude that
\begin{eqnarray}
\label{EqnLowerConc}
\mprob \big[ \sup_{v \in \SuperSpecSet{\qpar}} \big( 1 -\|X v\|_2
\big) \geq \frac{3}{2} \tlower(\newrad) \big] & \leq & 2 \exp
\big(-\numobs \frac{\tlower^2(\newrad)}{8} \big).
\end{eqnarray}

\newcommand{\TailTil}{\ensuremath{\widetilde{\Tail}}}

Define the event
\begin{eqnarray*}
\TailTil & \defn & \big \{ \exists \; v \in \real^\pdim \;
\mbox{s.t. $\|\CovMatSqrt v\|_2 = 1$ and } \big(1-\|X v \|_2) \geq 3
\tlower(\|v\|_1) \big \}.
\end{eqnarray*}
We can now apply Lemma~\ref{LemPeel} with $a_\numobs = \numobs$,
$g(\newrad) = 3 \tlower(\newrad)/2$ and $\mycon = 1/2$ to conclude
that there exist constants $\plaincon_i$ such that $\mprob[\TailTil]
\leq \plaincon_1 \exp(-\plaincon_2 \numobs)$.

Finally, to extend the claim to all vectors $v$, we consider the
rescaled vector $\vbre = v/\|\CovMatSqrt v\|_2$.  Conditioned on the
event $\TailTil^c$, we have for all $v \in \real^d$,
\begin{eqnarray*}
1-\|X \vbre\|_2/\sqrt{\numobs} & \leq & \frac{3}{2}  + 3 \,
(\sqrt{(3 \rhomax(\CovMat) \,\log \pdim)/\numobs} )\; \|\vbre\|_1,
\end{eqnarray*}
or equivalently, after multiplying through by $\|\CovMatSqrt v\|_2$
and re-arranging,
\begin{eqnarray*}
\|X v\|_2/\sqrt{\numobs} & \geq & \frac{1}{2} \|\CovMatSqrt v\|_2
  - 3 \, (\sqrt{(3 \rhomax(\CovMat) \,\log \pdim)/\numobs} )\; \|v\|_1,
\end{eqnarray*}
as claimed.
%

%%%%%%%%%%%%%%%%%%%%%%%%%%%%%%%%%%%%%%%%%%%%%%%%%%%%%%%%%%%%%%%%%%%%%%%%%%%%%%

%%%%%%%%%%%%%%%%%%%%%%%%%%%%%%%%%%%%%%%%%%%%%%%%%%%%%%%%%%%%%%%%%%%%%%%%%%%

\section{Proof of Lemma 6}
\label{AppLemNewrad}

For a given radius $\radtwo > 0$, define the set 
\begin{eqnarray*}
\InterSet(\spindex, \radtwo) & \defn & \big \{ \mytheta \in
\real^\pdim \, \mid \, \|\mytheta\|_0 \leq 2 \spindex, \quad
\|\mytheta\|_2 \leq \radtwo \big \},
\end{eqnarray*}
 and the random variables $Z_\numobs = Z_\numobs(\spindex, \radtwo)$
given by
\begin{eqnarray*}
Z_\numobs & \defn & \sup_{\theta \in \InterSet(\spindex, \radtwo)}
 \frac{1}{\numobs} |\wnoise^T X \theta|.
\end{eqnarray*}
For a given $\epsilon \in (0,1)$ to be chosen, let us upper bound the
minimal cardinality of a set that covers $\InterSet(\spindex,
\radtwo)$ up to $(\newrad \epsilon)$-accuracy in $\ell_2$-norm. We
claim that we may find such a covering set $\{ \theta^1, \ldots,
\theta^N\} \subset \InterSet(\spindex, \radtwo)$ with cardinality $N =
N(\spindex, \radtwo, \epsilon)$ that is upper bounded as
\begin{eqnarray*}
\log N(\spindex, \newrad, \epsilon) & \leq & \log {\pdim \choose 2
\spindex} + 2 \spindex \log(1/\epsilon).
\end{eqnarray*}
To establish this claim, note that here are ${\pdim \choose 2
\spindex}$ subsets of size $2 \spindex$ within $\{1, 2, \ldots,
\pdim\}$. Moreover, for any $2 \spindex$-sized subset, there is an
$(\radtwo \epsilon)$-covering in $\ell_2$-norm of the ball
$\Ball_2(\radtwo)$ with at most $2^{2 \spindex \log(1/\epsilon)}$
elements (e.g.,~\cite{Matousek}).

Consequently, for each $\theta \in \InterSet(\spindex, \radtwo)$, we
may find some $\theta^k$ such that $\|\theta - \theta^k\|_2 \leq
\radtwo \epsilon$.  By triangle inequality, we then have
\begin{eqnarray*}
\frac{1}{\numobs} |\wnoise^T X \theta| & \leq & \frac{1}{\numobs}
|\wnoise^T X \theta^k| + \frac{1}{\numobs} |\wnoise^T X (\theta -
\theta^i)| \\
& \leq & \frac{1}{\numobs} |\wnoise^T X \theta^k| +
\frac{\|\wnoise\|_2}{\sqrt{\numobs}} \, \frac{\|\Xmat (\theta -
\theta^k)\|_2}{\sqrt{\numobs}}.
\end{eqnarray*}
Given the assumptions on $\Xmat$, we have $\|\Xmat (\theta -
\theta^k)\|_2/\sqrt{\numobs} \leq \conupper \radtwo \|\theta - \theta^k\|_2 \; \leq \; \conupper \, \epsilon$.  Moreover, since the variate
$\|\wnoise\|_2^2/\sigma^2$ is $\chi^2$ with $\numobs$ degrees of
freedom, we have $\frac{\|\wnoise\|_2}{\sqrt{\numobs}} \leq 2 \sigma$
with probability $1 - \gencon_1 \exp(-\gencon_2 \numobs)$, using
standard tail bounds (see Appendix~\ref{AppChiTail}).  Putting
together the pieces, we conclude that
\begin{eqnarray*}
\frac{1}{\numobs} |\wnoise^T X \theta| & \leq & \frac{1}{\numobs}
|\wnoise^T X \theta^k| + 2 \conupper \, \sigma \, \radtwo \,
\epsilon
\end{eqnarray*}
with high probability. Taking the supremum over $\theta$ on both
sides yields
\begin{eqnarray*}
Z_\numobs & \leq & \max_{k = 1, 2, \ldots, N } \; \frac{1}{\numobs}
|\wnoise^T X \theta^k| + 2 \conupper \, \sigma \, \radtwo \,
\epsilon.
\end{eqnarray*}

It remains to bound the finite maximum over the covering set.  We
begin by observing that each variate $w^T \Xmat \theta^k/\numobs$ is
zero-mean Gaussian with variance $\sigma^2 \|\Xmat
\theta^i\|_2^2/\numobs^2$.  Under the given conditions on $\theta^k$
and $\Xmat$, this variance is at most $\sigma^2 \conuppersq
\radtwo^2/\numobs$, so that by standard Gaussian tail bounds, we
conclude that
\begin{eqnarray}
Z_\numobs & \leq & \sigma \; \radtwo \;  \conupper \, \sqrt{\frac{3 \log
N(\spindex, \radtwo, \epsilon)
}{\numobs}} + 2 \conupper \, \sigma \radtwo \, \epsilon \nonumber \\
\label{EqnMary}
& = & \sigma \; \radtwo \; \conupper \, \Big \{ \sqrt{\frac{3
\log N(\spindex, \radtwo, \epsilon)}{\numobs}} + 2 \epsilon \Big \}.
\end{eqnarray}
with probability greater than $1 - \gencon_1 \exp(-\gencon_2 \log
N(\spindex, \radtwo, \epsilon))$.

Finally, suppose that $\epsilon = \sqrt{\frac{\spindex \log(\pdim/2
\spindex)}{\numobs}}$.  With this choice and recalling that $\numobs
\leq \pdim$ by assumption, we obtain
\begin{eqnarray*}
\frac{\log N(\spindex, \radtwo, \epsilon)}{\numobs} & \leq &
\frac{\log {\pdim \choose 2 \spindex}}{\numobs} + \frac{\spindex \log
  \frac{\numobs}{\spindex \log(\pdim/2 \spindex)}}{\numobs} \\
& \leq & \frac{\log {\pdim \choose 2 \spindex}}{\numobs} +
  \frac{\spindex \log (\pdim/\spindex)}{\numobs} \\
& \leq & \frac{2 \spindex + 2 \spindex \log (\pdim/\spindex)} {\numobs}
  + \frac{\spindex \log (\pdim/\spindex)}{\numobs},
\end{eqnarray*}
where the final line uses standard bounds on binomial coefficients.
Since $\pdim/\spindex \geq 2$ by assumption, we conclude that our
choice of $\epsilon$ guarantees that $\frac{\log N(\spindex, \radtwo,
\epsilon)}{\numobs} \, \leq \, 5 \, \spindex \log(\pdim/\spindex)$.
Substituting these relations into the inequality~\eqref{EqnMary}, we
conclude that
\begin{eqnarray*}
Z_\numobs & \leq & \sigma \; \radtwo \; \conupper \, \Big \{ 4
\sqrt{\frac{\spindex \log(\pdim/\spindex)}{\numobs}} + 2
\sqrt{\frac{\spindex \log(\pdim/\spindex)}{\numobs}} \Big \},
\end{eqnarray*}
as claimed.  Since $\log N(\spindex, \radtwo, \epsilon) \geq \spindex
\log(\pdim - 2 \spindex)$, this event occurs with probability at least
$1 - \gencon_1 \exp(-\gencon_2 \min \{\numobs, \spindex
\log(\pdim-\spindex) \})$, as claimed.  

%%%%%%%%%%%%%%%%%%%%%%%%%%%%%%%%%%%%%%%%%%%%%%%%%%%%%%%%%%%%%%%%%%%%%%%%%%

\section{Proofs for Theorem 4}

This appendix is devoted to the proofs of technical lemmas
used in Theorem~\ref{ThmUpperPred}.

\subsection{Proof of Lemma~\ref{LemNewradPred}}
\label{AppLemNewradPred}

For $\qpar \in (0,1)$, let us define the set
\begin{eqnarray*}
\InterSet_\qpar(\myrad, \radtwo) & \defn & \Ballq(2 \myrad) \cap \big
\{ \mytheta \in \real^\pdim \, \mid \, \|\Xtil
\mytheta\|_2/\sqrt{\numobs} \leq \radtwo \big \}.
\end{eqnarray*}
We seek to bound the random variable $Z(\myrad, \radtwo) \defn
\sup_{\mytheta \in \InterSet_\qpar(\myrad, \radtwo)} \frac{1}{\numobs}
|\wtil^T \Xtil \mytheta|$, which we do by a chaining result---in
particular, Lemma 3.2 in van de Geer~\cite{vandeGeer}).  Adopting the
notation from this lemma, we seek to apply it with $\epsilon
=\delta/2$, and $K=4$.  Suppose that $\frac{\|X \theta\|_2}{\sqrt{n}}
\leq \radtwo$, and
\begin{subequations}
\begin{eqnarray}
\label{EqnVDGOne}
\sqrt{\numobs} \delta & \geq & \gencon_1 \radtwo \\ 
\label{EqnVDGTwo}
\sqrt{\numobs} \delta & \geq & \gencon_1
\int_{\frac{\delta}{16}}^\radtwo \sqrt{\log
\CovNum(t;\InterSet_\qpar)} dt \, = \, : J(\radtwo, \delta).
\end{eqnarray}
\end{subequations}	 
where $\CovNum(t;\InterSet_\qpar)$ is the covering number for
$\InterSet_\qpar$ in the $\ell_2$-prediction norm (defined by $\|\Xmat
\theta\|/ \sqrt{\numobs}$).  As long as
\mbox{$\frac{\|\wtil\|_2^2}{\numobs} \leq 16$,} Lemma 3.2 guarantees
that
\begin{eqnarray*}
\mprob\big[Z(\myrad, \radtwo) \geq \delta, \;
\frac{\|\wtil\|_2^2}{\numobs} \leq 16] \leq \gencon_1 \exp{(-\gencon_2
\frac{\numobs \delta^2}{\radtwo^2})} .
\end{eqnarray*}
By tail bounds on $\chi^2$ random variables (see
Appendix~\ref{AppChiTail}), we have $\mprob[\|\wtil\|_2^2 \geq 16
\numobs] \leq \gencon_4 \exp(-\gencon_5 \numobs)$.  
Consequently, we conclude that
\begin{align*}
\mprob\big[Z(\myrad, \radtwo) \geq \delta] & \leq \gencon_1
\exp{(-\gencon_2 \frac{\numobs \delta^2}{\radtwo^2})} + \gencon_4
\exp(-\gencon_5 \numobs)
\end{align*}

For some $\gencon_3 > 0$, let us set
\begin{eqnarray*}
\delta & = & \gencon_3 \: \radtwo \: \colnormtil^{\frac{\qpar}{2}} \,
\sqrt{\myrad} \; (\frac{\log \pdim}{\numobs})^{\frac{1}{2} -
\frac{\qpar}{4}},
\end{eqnarray*}
and let us verify that the conditions~\eqref{EqnVDGOne}
and~\eqref{EqnVDGTwo} hold.  Given our choice of $\delta$, we find
that
\begin{equation*}
\frac{\delta}{\radtwo} \, \sqrt{\numobs} = \Omega(\numobs^{\qpar/4}
(\log \pdim)^{1/2 - \qpar/4}),
\end{equation*}
and since $\pdim, \numobs \rightarrow \infty$, we see that
condition~\eqref{EqnVDGOne} holds.  Turning to verification of the
inequality~\eqref{EqnVDGTwo}, we first provide an upper bound for
$\log N(\InterSet_\qpar, t)$. Setting $\gamma = \frac{\Xtil
\mytheta}{\sqrt{\numobs}}$ and from the definition~\eqref{EqnQparConv}
of $\absconv_\qpar(\Xmat/\sqrt{\numobs})$, we have
\begin{equation*} 
\sup_{\mytheta \in \InterSet_\qpar(\myrad, \radtwo)} \frac{1}{\numobs}
|\wtil^T \Xtil \mytheta| \leq \sup_{\gamma \in
  \absconv_\qpar(\Xmat/\sqrt{\numobs}), \|\gamma\|_2 \leq \radtwo}
\frac{1}{\sqrt{\numobs}} |\wtil^T \gamma|.
\end{equation*}
We may apply the bound in Lemma~\ref{LemKL} to conclude that $\log
\CovNum(\epsilon;\InterSet_\qpar)$ is upper bounded by $\gencon' \;
\Rq^{\frac{2}{2-\qpar}} \; \big(\frac{\colnormtil}{\epsilon}
\big)^{\frac{2 \qpar}{2-\qpar}} \log \pdim$.

Using this upper bound, we have
\begin{eqnarray*}
J(\radtwo, \delta) \defn \int_{\delta/16}^\radtwo \sqrt{\log
N(\InterSet_\qpar, t)} dt & \leq & \int_{0}^\radtwo \sqrt{\log
N(\InterSet_\qpar, t)} dt\\ & \leq & \gencon \, \, \;
\Rq^{\frac{1}{2-\qpar}} \; \colnormtil^{\frac{\qpar}{2-\qpar}} \;
\sqrt{\log \pdim} \int_{0}^\radtwo t^{-\qpar/(2-\qpar)} dt \\
& = & \gencon' \Rq^{\frac{1}{2-\qpar}} \;
\colnormtil^{\frac{\qpar}{2-\qpar}} \; \sqrt{\log \pdim} 
\; \radtwo^{1-\frac{\qpar}{2-\qpar}}.
\end{eqnarray*}

Using this upper bound, let us verify that the
inequality~\eqref{EqnVDGTwo} holds as long as $\radtwo =
\Omega(\colnormtil^{\frac{\qpar}{2}} \, \sqrt{\myrad} \; (\frac{\log
\pdim}{\numobs})^{\frac{1}{2} - \frac{\qpar}{4}})$, as assumed in the
statement of Lemma~\ref{LemNewradPred}.  With our choice of $\delta$,
we have
\begin{align*}
\frac{J}{\sqrt{\numobs} \, \delta} & \leq 
\frac{\gencon'
\Rq^{\frac{1}{2-\qpar}} \; \colnormtil^{\frac{\qpar}{2-\qpar}} \;
\sqrt{\frac{\log \pdim}{\numobs}} \;
\radtwo^{1-\frac{\qpar}{2-\qpar}}}{ \gencon_3 \: \radtwo \:
\colnormtil^{\frac{\qpar}{2}} \, \sqrt{\myrad} \; (\frac{\log
\pdim}{\numobs})^{\frac{1}{2} - \frac{\qpar}{4}}} \\
& = \frac{\gencon' \Rq^{\frac{1}{2-\qpar} - \frac{1}{2} -
\frac{\qpar}{2 \,(2-\qpar)}} \; \colnormtil^{\frac{\qpar}{2-\qpar} -
\frac{\qpar}{2} \frac{\qpar}{2-\qpar} - \frac{\qpar}{2} } \;
\big(\frac{\log \pdim}{\numobs}\big)^{\frac{\qpar}{4} -
\frac{\qpar}{2-\qpar} \big(\frac{1}{2} - \frac{\qpar}{4} \big)}} {
\gencon_3} \\
& = \frac{\gencon'}{\gencon_3},
\end{align*}
so that condition~\eqref{EqnVDGTwo} will hold as long as we choose
$\gencon_3 > 0$ large enough.  Overall, we conclude that $\mprob[
Z(\myrad, \radtwo) \geq \gencon_3 \: \radtwo \:
\colnormtil^{\frac{\qpar}{2}} \, \sqrt{\myrad} \; (\frac{\log
\pdim}{\numobs})^{\frac{1}{2} - \frac{\qpar}{4}}] \leq \gencon_1
\exp(-\Rq (\log \pdim)^{1-\frac{\qpar}{2}} n^{\frac{\qpar}{2}})$,
which concludes the proof.

%%%%%%%%%%%%%%%%%%%%%%%%%%%%%%%%%%%%%%%%%%%%%%%%%%%%%%%%%%%%%%%%%%%%%%%

\subsection{Proof of Lemma~\ref{LemNewradPredZero}}
\label{AppLemNewradPredZero}
First, consider a fixed subset $\Sset \subset \{1, 2, \ldots, \pdim
\}$ of cardinality $|\Sset| = \spindex$.  Applying the SVD to the
sub-matrix $\Xmat_\Sset \in \real^{\numobs \times \spindex}$, we have
$\Xmat_\Sset = V D U$, where $V \in \real^{\numobs \times \spindex}$
has orthonormal columns, and $D U \in \real^{\spindex \times
\spindex}$.  By construction, for any $\Delta_\Sset \in
\real^\spindex$, we have $\|\Xmat_\Sset \Delta_\Sset\|_2 = \|D U
\Delta_\Sset\|_2$.  Since $V$ has orthonormal columns, the vector
$\wtil_\Sset = V^T w \in \real^\spindex$ has i.i.d. $N(0, \sigma^2)$
entries.  Consequently, for any $\Delta_\Sset$ such that
$\frac{\|\Xmat_\Sset \Delta_\Sset\|_2}{\sqrt{\numobs}} \leq \radtwo$,
we have
\begin{eqnarray*}
\Big| \frac{w^T \Xmat_\Sset \Delta_\Sset}{\numobs} \Big | & = & \Big |
\frac{\wtil^T_\Sset}{\sqrt{\numobs}} \frac{D U
\Delta_\Sset}{\sqrt{\numobs}} \Big| \\
& \leq & \frac{\|\wtil_\Sset\|_2}{\sqrt{\numobs}} \frac{\|D U
\Delta_\Sset \|_2}{\sqrt{\numobs}} \\
& \leq & \frac{\|\wtil_\Sset\|_2}{\sqrt{\numobs}} \; \radtwo.
\end{eqnarray*}
Now the variate $\sigma^{-2} \|\wtil_\Sset\|_2^2$ is $\chi^2$ with
$\spindex$ degrees of freedom, so that by standard $\chi^2$ tail
bounds (see Appendix~\ref{AppChiTail}), we have
\begin{eqnarray*}
\mprob\big[\frac{\|\wtil_\Sset\|_2^2}{\sigma^2 \spindex} \geq 1 + 4
\delta \big] & \leq & \exp(-\spindex \delta), \mbox{ valid for all
$\delta \geq 1$.}
\end{eqnarray*}
 Setting $\delta = 20 \log(\frac{\pdim}{2 \spindex})$ and noting that
$\log(\frac{\pdim}{2 \spindex}) \geq \log 2$ by assumption, we have
(after some algebra)
\begin{eqnarray*}
\mprob \biggr[\frac{\|\wtil_\Sset\|_2^2}{\numobs} \geq \frac{\sigma^2
\spindex}{\numobs} \big(81 \log (\pdim/\spindex) \big) \biggr] &
\leq & \exp(-20 \spindex \log(\frac{\pdim}{2 \spindex})).
\end{eqnarray*}
We have thus shown that for each fixed subset, we have the bound
\begin{eqnarray*}
\Big |\frac{w^T \Xmat_\Sset \Delta_\Sset}{\numobs} \Big | & \leq &
\radtwo \, \sqrt{\frac{81 \sigma^2 \spindex \log(\frac{\pdim}{2
\spindex})}{\numobs}},
\end{eqnarray*}
with probability at least $1 - \exp(-20 \spindex \log(\frac{\pdim}{2
\spindex}))$.

Since there are ${\pdim \choose 2 \spindex} \leq
\big(\frac{\pdim e}{2 \spindex})^{2 \spindex}$ subsets of size
$\spindex$, applying a union bound yields that
\begin{eqnarray*}
\mprob \Big[ \sup_{\mytheta \in \Ball_0(2 \spindex), \; \frac{\|\Xmat
\mytheta\|_2}{\sqrt{\numobs}} \leq \radtwo} |\frac{w^T \Xmat
\mytheta}{\numobs}| \geq \radtwo \, \sqrt{\frac{81 \sigma^2 \spindex
\log(\frac{\pdim}{2 \spindex})}{\numobs}} \Big] & \leq & \exp
\biggr( -20 \spindex \log(\frac{\pdim}{2 \spindex}) + 2 \spindex \log
\frac{\pdim e}{2 \spindex} \biggr) \\
& \leq & \exp \big( -10 \spindex \log(\frac{\pdim}{2 \spindex}) \big),
\end{eqnarray*}
as claimed.

\section{Large deviations for random objectives}
\label{AppLemPeel}

In this appendix, we state a result on large deviations of the
constrained optimum of random objective functions of the form $f(v;
X)$, where $v \in \real^\pdim$ is the optimization vector, and $X$ is
some random vector.  Of interest is the optimization problem
$\sup_{\rho(v) \leq \newrad, \; v \in S} f(v; X_n)$, where $\rho:
\real^\pdim \rightarrow \real_+$ is some non-negative and increasing
constraint function, and $S$ is a non-empty set.  With this set-up,
our goal is to bound the probability of the event defined by
\begin{eqnarray*}
\Event & \defn & \big \{ \exists \; v \in S \mbox{ such that } f(v; X)
 \geq 2 g(\rho(v))) \big \},
\end{eqnarray*}
where $g: \real \rightarrow \real$ is non-negative and strictly
increasing. 
\blems
\label{LemPeel}
Suppose that $g(\newrad) \geq \mycon$ for all $\newrad \geq 0$, and
that there exists some constant $\gencon > 0$ such that for all
$\newrad > 0$, we have the tail bound
\begin{eqnarray}
\label{EqnBasicTail}
\mprob \big[ \sup_{v \in S, \; \rho(v) \leq \newrad} f_n(v; X_n) \geq
g(\newrad)] & \leq & 2 \exp(-\gencon \: a_n \: g^2(\newrad)),
\end{eqnarray}
for some $a_n > 0$.  Then we have
\begin{eqnarray}
\mprob[\Event_n] & \leq & \frac{2 \exp(-\gencon a_n \mycon^2)}{1 -
\exp(- \gencon a_n \mycon^2)}.
\end{eqnarray}
\elems
\spro
Our proof is based on a standard peeling technique (e.g., see van de
Geer~\cite{vandeGeer} pp.~82).  By assumption, as $v$ varies over $S$,
we have $g(\newrad) \in [\mycon, \infty)$.  Accordingly, for $m = 1,
2, \ldots$, defining the sets
\begin{eqnarray*}
S_m & \defn & \big \{ v \in S \, \mid \, 2^{m-1} \mycon \leq
g(\rho(v)) \leq 2^m \mycon \big \},
\end{eqnarray*}
we may conclude that if there exists $v \in S$ such that $f(v, X) \geq
2 h(\rho(v))$, then this must occur for some $m$ and $v \in S_m$.  By
union bound, we have
\begin{eqnarray*}
\mprob[\Event] & \leq & \sum_{m=1}^\infty \mprob \big[\exists \; v \in
S_m \mbox{ such that } f(v, X) \geq 2 g(\rho(v)) \big].
\end{eqnarray*}
If $v \in S_m$ and $f(v, X) \geq 2 g(\rho(v))$, then by definition of
$S_m$, we have \mbox{$f(v, X) \geq 2 \, (2^{m-1}) \, \mycon = 2^m
\mycon$.}  Since for any $v \in S_m$, we have $g(\rho(v)) \leq 2^m
\mycon$, we combine these inequalities to obtain
\begin{eqnarray*}
\mprob[\Event] & \leq & \sum_{m=1}^\infty \mprob \big[\sup_{\rho(v)
\leq g^{-1}(2^m \mycon)} f(v, X) \geq 2^m \mycon \big] \\
& \leq & \sum_{m=1}^\infty 2 \exp\big(-\gencon a_\numobs \; [g(g^{-1}(2^m
\mycon))]^2 \big) \\
& = & 2 \sum_{m=1}^\infty \exp\big(-\gencon a_\numobs \; 2^{2m}
 \mycon^2 \big),
\end{eqnarray*}
from which the stated claim follows by upper bounding this geometric
sum.
\fpro

%%%%%%%%%%%%%%%%%%%%%%%%%%%%%%%%%%%%%%%%%%%%%%%%%%%%%%%%%%%%%%%%%%%%%%%%%%

\section{Some tail bounds for $\chi^2$-variates}
\label{AppChiTail}

\newcommand{\dfchi}{\ensuremath{m}}

The following large-deviations bounds for centralized $\chi^2$ are
taken from Laurent and Massart~\cite{LauMas98}.  Given a centralized
$\chi^2$-variate $Z$ with $\dfchi$ degrees of freedom, then for all $x
\geq 0$,
\begin{subequations}
\begin{eqnarray}
\label{EqnCleanUpCent}
\mprob \left[Z - \dfchi \geq 2 \sqrt{\dfchi x} + 2x \right] & \leq &
\exp(-x), \qquad \mbox{and} \\
\label{EqnCleanDownCent}
\mprob \left[Z - \dfchi \leq -2 \sqrt{\dfchi x} \right] & \leq &
\exp(-x).
\end{eqnarray}
\end{subequations}
The following consequence of this bound is useful: for $t \geq 1$, we
have
\begin{eqnarray}
\label{EqnChitailBig}
\mprob \big[\frac{Z - \dfchi}{\dfchi} \geq 4 t \big] & \leq &
\exp(-\dfchi t).
\end{eqnarray}
Starting with the bound~\eqref{EqnCleanUpCent}, setting $x = t \dfchi$
yields $\mprob \big[\frac{Z - \dfchi}{\dfchi} \geq 2 \sqrt{t} + 2 t
\big] \, \leq \, \exp(-t \dfchi)$, Since $4 t \geq 2 \sqrt{t} + 2t$
for $t \geq 1$, we have $\mprob[\frac{Z - \dfchi}{\dfchi} \geq 4 t]
\leq \exp(-t \dfchi)$ for all $t \geq 1$. 

%%%%%%%%%%%%
%\bibliographystyle{plainnat}
%\bibliography{GarveshAug09}

%%%%%%%%%%%%%%%%

\end{document}